\documentclass{article}[11pt]
\usepackage[a4paper,top=4.0cm,bottom=2.0cm,left=2.0cm,right=2.0cm]{geometry}

\usepackage{amsmath,amsfonts,mathrsfs,amsfonts,ulem}
\usepackage{epsfig,graphicx}
\usepackage{color}

\usepackage{amsmath,amsfonts,mathrsfs,amsfonts, color, bbold,bm,blindtext}

\usepackage{algorithmicx,algorithm,algpseudocode}

\allowdisplaybreaks
\begin{document}
	
	\title{Ensemble transform algorithms for nonlinear smoothing problems}

	\author{Jana de Wiljes\thanks{Universit\"at Potsdam, 
			Institut f\"ur Mathematik, Karl-Liebknecht-Str. 24/25, D-14476 Potsdam, Germany ({\tt wiljes@uni-potsdam.de})  and University of Reading, Department of Mathematics and Statistics, Whiteknights, PO Box 220, Reading RG6 6AX, UK} \and 
		Sahani Pathiraja\thanks{Universit\"at Potsdam, 
			Institut f\"ur Mathematik, Karl-Liebknecht-Str. 24/25, D-14476 Potsdam, Germany ({\tt pathiraja@uni-potsdam.de})}\and
		Sebastian Reich\thanks{Universit\"at Potsdam, Institut f\"ur Mathematik, Karl-Liebknecht-Str. 24/25, D-14476 Potsdam, Germany ({\tt sreich@math.uni-potsdam.de}) and University of Reading, Department of Mathematics and Statistics, Whiteknights, PO Box 220, Reading RG6 6AX, UK} 
	}

	\maketitle
	
	\begin{abstract} 
	
		Several numerical tools designed to overcome the challenges of smoothing in a nonlinear and non-Gaussian setting are investigated for a class of particle smoothers. The considered family of smoothers is induced by the class of linear ensemble transform filters  which contains classical filters such as the stochastic ensemble Kalman filter, the ensemble square root filter and the recently introduced nonlinear ensemble transform filter. Further the ensemble transform particle smoother is introduced and particularly highlighted as it is consistent in the particle limit and does not require assumptions with respect to the family of the posterior distribution. The linear update pattern of the considered class of linear ensemble transform smoothers allows one to implement important supplementary techniques such as adaptive spread corrections, hybrid formulations, and localization in order to facilitate their application to complex estimation problems. These additional features are derived and numerically investigated for a sequence of increasingly challenging test problems.
		
	\end{abstract}
	
	\noindent
	{\bf Keywords.} Data assimilation, particle smoother, optimal transport, localisation, hybrid, second order correction\\
	\noindent {\bf AMS(MOS) subject classifications.} 65C05, 62M20, 93E11, 62F15, 86A22

	%%%%%%%%%%%%%%%%%%%%%%%%%%%%%%%%%%%%%%%%%
	%
	\section{Introduction}
%%%%%%%%%%%%%%%%%%%%%%%%%%%
	In many applicational areas, e.g., in meteorology \cite{jdw:Shuman1957}, neuroscience \cite{Hu2018} and robotics \cite{saerkkae2013}, there is a need to smooth model states using data extending further into the future. The underlying setting of such state estimation problems consists of two information sources: a model describing the system of interest and observations of the model state. Here we consider a discrete-time Markov state space model of the general form
	\begin{equation} \label{eq:t1}
	\begin{aligned}
	{\bf x}_{k} = f({\bf x}_{k-1}) + \epsilon_k, \ \epsilon_k \sim \mathcal{N}({\bf 0},{\bf Q}) 
	\end{aligned}
	\end{equation}
	where $ {\bf x}_{k} \in \mathbb{R}^{N_x} $ denotes the state of the system at time $k$, $ {\bf Q} \in \mathbb{R}^{N_x \times N_x}$ is the positive semi-definite covariance matrix of model errors and the initial states ${\bf x}_{0}$ are drawn from a given distribution $p({\bf x}_{0})$. Here 
	$\mathcal{N}({\bf 0},{\bf Q})$ stands for the Gaussian distribution with mean zero and covariance matrix ${\bf Q}$ and we also allow for the case ${\bf Q} = {\bf 0}$
	in which case the model (\ref{eq:t1}) becomes deterministic. Additionally, partial and noisy measurements of the system states are available according to the following process:
	\begin{equation} \label{eq:t2}
	\begin{aligned}
	{\bf y}_{k} = h({\bf x}_{k}) + \nu_k,  \ \nu_k \sim \mathcal{N}({\bf 0},{\bf R}) 
	\end{aligned}
	\end{equation}
	where ${\bf y}_{k} \in \mathbb{R}^{N_y}$ is the measurement at time $k$ for $k=1,\ldots,K$, and ${\bf R} \in \mathbb{R}^{N_y \times N_y}$ is a positive definite matrix of measurement errors. The notation
	\begin{equation} 
	{\bf y}_{1:k} = [{\bf y}_{1}^{\rm T}, {\bf y}_{2}^{\rm T}, ..., {\bf y}_{k}^{\rm T}]^{\rm T} \in \mathbb{R}^{k N_y \times 1}
	\end{equation}
	 is used to indicate the collection of observations from time one to $k \le K$. The same short-hand notation will be used for collection of state vectors.
	The overall goal in the resulting estimation problems is usually to either approximate the marginal filtering $p({\bf x}_{k}|{\bf y}_{1:k})$ or smoothing 
	distributions $p({\bf x}_l|{\bf y}_{1:k})$ for $l = 0,\ldots,k-1$.\footnote{With a slight abuse of notation we use the same symbol $p$ to denote different distributions. It should be clear from the arguments in $p$ which specific distribution is referred to.} Whilst filtering is useful for real-time prediction where the data comes in sequentially, smoothers are more desirable when one has the possibility to asses a whole time interval of measurements at once and can therefore obtain more informative posteriors by considering the full data set. This is particularly effective for state estimation in the context of models with long memory or long temporal correlations, e.g.~ocean circulation modeling \cite{jdw:Cosmeetal2009,jdw:Ocean2006}.
In some applications, one might also be interested in the joint smoothing distribution 
	$p({\bf x}_{k-L:k}|{\bf y}_{1:k})$ over some time window of length $L\le k$. This is, for example, the case when the discrete-time Markov model
	(\ref{eq:t1}) is replaced by a model of order $n>1$ of the general form
	\begin{equation} \label{eq:t1b}
	\begin{aligned}
	{\bf x}_{k} = f({\bf x}_{k-n},\ldots,{\bf x}_{k-1}) + \epsilon_k, \ \epsilon_k \sim \mathcal{N}({\bf 0},{\bf Q}) ,
	\end{aligned}
	\end{equation}
	in which case sequential state estimation requires knowledge of the joint smoothing distribution 
	$p({\bf x}_{k-n:k-1}|{\bf y}_{1:k-1})$ at time $k-1$, that is $L=n$, in order to proceed to time $k$
	using the model (\ref{eq:t1b}) and the data (\ref{eq:t2}).

Due to these advantages, smoothing is of fundamental importance in several applicational areas. For example variational approaches \cite{jdw:FisherLeutbecherKelly2005} such as 4DVAR \cite{jdw:LeDimetTalagrand1986} are well established in Numerical Weather Prediction and randomized maximum likelihood and iterative ensemble Kalman smoothers \cite{Bocquet2014,jdw:ChenOliver2017,jdw:Evensen2018} are key methods for oil reservoir modelling. Our focus is strictly on derivative-free ensemble based smoothing techniques as opposed to the large family of methods that rely on minimisation of certain cost functionals such as the aforementioned methods.

In the linear-Gaussian case, the smoothing distribution remains Gaussian with mean and covariance that can be computed by the Kalman smoother, which is a direct generalization of the Kalman filter \cite{jdw:Kalman1960,sr:jazwinski}. 
In the nonlinear setting, the popular family of ensemble Kalman filters (EnKFs) \cite{sr:evensen2006} has inspired a number of different smoothers \cite{jdw:EvensenLeeuwen2000,jdw:Cosmeetal2009,jdw:WhitakerCompo2002}. The reader is referred to \cite{jdw:Carrassietall2018review} and references therein for a recent review on the various existing approaches.

Alternatively there are several smoothing extensions of particle filters \cite{saerkkae2013}. The canonical ansatz for a particle smoother is to build on the classical bootstrap filter to obtain smoothed state space estimates with respect to observations extending beyond the current time \cite{jdw:DoucetGodsillAndrieu2000,jdw:KitagawaSato2001}.  However particle degeneracy is even more prominent in the smoothing case than in filtering. If the smoother is only applied over a relatively small time interval, degeneracy can be minimised significantly \cite{jdw:ClappGodsill1998}. Further, a number of different implementations, e.g., using variance reduction techniques, have been proposed and a nice overview can be found in \cite{sr:Doucet,saerkkae2013}. Further solutions to the smoothing problem can also be derived using optimal control techniques by assuming a linear evolution model for the trajectories \cite{jdw:DoucetGodsillAndrieu2000,sr:Doucet}.
Moreover smoothing approaches have recently been used to improve popular filtering techniques such as the auxiliary particle filter \cite{jdw:PetetinDesbouvries2011} or the implicit particle method in \cite{Weiretal2013}. Conceptually the idea is to utilize proposal densities that factor in observations at future time steps (typically one data point beyond the current time instance is already sufficient to improve the underlying filter).
%Smoothing can be viewed as a form of filtering with an augmented state space (where the vector of states at the current time is augmented with states in the past) which increases the dimensionality of the problem.

Here we propose a unifying framework, referred to as linear ensemble transform smoothers (LETSs), to highlight the connection between various smoothing algorithms.  As we will show, the family of LETSs includes a range of popular smoothers such as the class of ensemble square root smoothers (ESRS) \cite{jdw:EvensenLeeuwen2000, Cosme2010, jdw:WhitakerCompo2002} and the nonlinear ensemble transform smoother (NETS) \cite{sr:toedter15}.  We propose for the first time, a smoother extension of the ensemble transform particle filter (ETPF) \cite{sr:reichcotter15}, which requires a non-trivial modification of existing LETS formulations.  The advantage of the ensemble transform particle smoother (ETPS) is that it has no underlying Gaussianity assumptions and is formally consistent with the true posterior in the ensemble limit ({\it i.e.}, $M\rightarrow\infty$). 

The other main contribution of this manuscript is to present a series of tools that can be used to address computational difficulties associated with smoothing in moderate dimensional nonlinear problems, particularly for small to moderate ensemble sizes.  This is motivated by the fact that in many applications, particularly numerical weather prediction, ensemble methods with very small ensemble sizes relative to the state dimension are utilised \cite{Bonavita2008, Miyoshi2011, Nurujjaman2013}. The first tool we discuss is the so-called hybrid formulation, which allows one to develop superior smoothers by combining very robust but inconsistent smoothers such as the ESRS with consistent smoothers such as the ETPS.  We demonstrate how such hybrid formulations can be developed in the context of LETSs, taking inspiration from hybrid formulations recently introduced in the context of filters in \cite{sr:CRR15}.  Another important practical issue is that one is typically forced to use a small ensemble in many applications. A common problem associated with this in filtering is that the posterior spread is often underestimated. In practice this issue is approached via a random increase of the ensemble variance. A more judicious strategy is to match the second moment of ensemble to a spread one trusts to be better than the current approximation. With this in mind, we extend a moment matching technique described in \cite{acevedo2016} to the class of smoothers in this paper. Localisation is another popular strategy to deal with effects of spurious correlations due to small ensemble sizes \cite{sr:hunt07, sr:houtekamer01}, and can also be seen as a way to handle the curse of dimensionality in the particle filtering context \cite{jdw:FarchiBocquet2018}. 
 This strategy has also been utilized in the context of smoothing, e.g., see recent work in \cite{Morzfeld2017, Bocquet2016, jdw:ChenOliver2017}. Whilst in filtering one is typically concerned with localising in space, in the more general smoothing problem it is pertinent to also consider localisation in time.  In this paper, we propose a number of spatio-temporal localisation techniques specifically for LETSs.  Finally, we demonstrate the benefits of these approaches through a range of numerical experiments.

The remainder of the paper is structured as follows. First an overview of the mathematical background and notation is given. In Section \ref{sec:LETF} 
we first summarize the key idea behind the ESRS and the NETS and demonstrate why their straightforward extension to the ETPF fails. We then introduce a generalised family of LETSs in the form of (\ref{eq:LETS_proper}) which overcomes this limitation. This is followed by a discussion on numerical tools such as localisation in time, ensemble spread correction and hybrid formulations that can be employed in the context of the proposed LETSs.
The presented class of smoothers is numerically investigated for a range of systems in Section \ref{sec:numerics}: in particular we consider a higher order Markov system, as the effects of smoothing are particularly prominent when the underlying model has long memory.  We also consider a classical moderate dimensional toy system, since finite ensemble sizes are particularly problematic for high dimensional applications. Lastly in Section \ref{sec:spatially_extended}, we discuss how spatio-temporal localisation techniques can be utilized to make the smoothing algorithms feasible in that setting. 

%%%%%%%%%%%%%%%%%%%%%%%%%%%%%%%%%%%%%%%%%%%%%%%%%%%%%%%%%%%%%%%
		\section{Mathematical background and notation}\label{sec:Mathematicalframework}
%%%%%%%%%%%%%%%%%%%%%%%%%%%%%%%%%%%%%%%%%%%%%%%%%%%%%%%%%%%%%%%%

As the smoothing algorithms we will consider throughout the paper are implemented in a sequential manner it is important to understand the following equalities: The joint smoothing distribution $p({\bf x}_{0:k}|{\bf y}_{1:k})$ can simply be linked to the previous smoothing distribution  $p({\bf x}_{0:k-1}|{\bf y}_{1:k-1})$ using the assumption that model and measurement errors in (\ref{eq:t1}) and (\ref{eq:t2}), respectively, are independent:
	\begin{equation} 
	\label{b1}
	p({\bf x}_{0:k}|{\bf y}_{1:k}) = \frac{p({\bf y}_k|{\bf x}_k) p({\bf x}_k|{\bf x}_{k-1})}{p({\bf y}_k|{\bf y}_{1:k-1})} p({\bf x}_{0:k-1}|{\bf y}_{1:k-1}) 
	= \left[\prod^k_{l=1} \frac{p({\bf y}_l|{\bf x}_{l})}{p({\bf y}_{l}|{\bf y}_{1:l-1})}\cdot p({\bf x}_{l}|{\bf x}_{l-1}) \right] p({\bf x}_{0}).
\end{equation}
Here $p({\bf x}_k|{\bf x}_{k-1})$ denotes the transition probability associated with the Markov model (\ref{eq:t1}) and 
\begin{equation} \label{eq:likelihood}
p({\bf y}_k|{\bf x}_k) \propto \exp\Big(-\frac{1}{2}(h({\bf x}_{k})-{\bf y}_{k})^{\rm T}{\bf R}^{-1}(h({\bf x}_{k})-{\bf y}_{k})\Big)
\end{equation}
denotes the likelihood function associated with the measurement model (\ref{eq:t2}). The joint smoothing distribution 
(\ref{b1}) gives rise to the associated marginal distributions $p({\bf x}_l|{\bf y}_{1:k})$ in the state variable ${\bf x}_l$ 
for $l=k$ (filtering), $l<k$ (smoothing), or the fixed-lag smoothing distribution $p({\bf x}_{k-L:k}|{\bf y}_{1:k})$. Computing expectation values with 
respect to the joint smoothing distribution (\ref{b1}) or any of its marginal distributions is generally intractable except for the linear Gaussian case.  
We therefore consider sample-based Monte Carlo methods\footnote{This paper being at the crossroad of the ensemble Kalman and the sequential 
Monte Carlo/particle filter/smoothing literature, we use the notion of particles, samples, and ensemble members synonymously throughout
this paper.} in this paper and use the following notation to indicate $M$ samples 
${\bf x}_{l|k}^{(i)} \in \mathbb{R}^{N_x}$, $i=1,\ldots,M$, representing the marginal prediction ($l>k$), 
filtering ($l=k$), and smoothing ($l<k$) distributions respectively, which we write with a slight abuse of notation as
		\begin{align}
		{\bf x}^{(i)}_{l|k}  &\sim p({\bf x}_l|{\bf y}_{1:k}) ,
		\end{align}
where $0\le l \le k+1$ for $k = 0,\ldots,K$. Further, the full ensemble of samples at any time $k$ is collected in an $N_x\times M$ matrix such as
\begin{equation}
 {\bf X}_{l|k}  = \left({\bf x}^{(1)}_{l|k}, ..., {\bf x}^{(M)}_{l|k} \right) .
 \end{equation} 
One can also collect past and current state values of the $i$th ensemble member in a column vector ${\bf x}_{k-L:k|k}^{(i)}$ of dimension
$(L+1) N_x$ for any $1\le L \le k$ in order to produce samples from the associated fixed-lag $L$ smoothing distributions, which we 
write again with a slight abuse of notation as
\begin{equation} 
{\bf x}_{k-L:k|k}^{(i)} \sim p({\bf x}_{k-L:k}|{\bf y}_{1:k}).
\end{equation}
Again the full ensemble of such trajectory samples can be collected in an $(L+1)N_x \times M$ matrix
\begin{equation} \label{eq:matrix_trajectories}
 {\bf X}_{k-L:k|k}  = \left({\bf x}^{(1)}_{k-L:k|k}, ..., {\bf x}^{(M)}_{k-L:k|k} \right) .
 \end{equation}

Below we describe a classical importance sampling implementation of the particle smoother \cite{jdw:GodsillDoucetWest2012}.  The crucial difference to filtering is that after every assimilation step, the entire history is updated according to the new importance weight. Just as in filtering, it is important to perform resampling to mitigate weight degeneracy, although in smoothing, resampling is typically undertaken trajectory-wise \cite{saerkkae2013}. 
A single cycle of a classical particle smoother at time $k$ is therefore structured as follows:\medskip\\
\noindent\textbf{Prediction:}\footnote{We mention that the prediction step can utilize alternative proposal densities \cite{Arulampalam2002}. However, we will
not discuss such alternative proposals in this paper as they add another layer of computational complexity while generally not being able to beat
the curse of dimensionality of a filtering/smoothing problem \cite{Snyder2008, Snyder2015, MorzfeldHodyssSnyder2017}.} Given $M$ samples ${\bf x}_{0:k-1|k-1}^{(i)} \in \mathbb{R}^{kN_x}$ 
from the smoothing distribution at time $k-1$, one first produces a 
forecast sample ${\bf x}_{k|k-1}^{(i)} \in \mathbb{R}^{N_x}$ 
at time $k$ using (\ref{eq:t1}) with ${\bf x}_{k-1} = {\bf x}_{k-1|k-1}^{(i)}$ for all $i=1,\ldots,M$. This leads to the concatenated trajectory
samples ${\bf x}^{(i)}_{0:k|k-1} \in \mathbb{R}^{(k+1)N_x}$ and the prior density for the subsequent Bayesian inference step is approximated 
via the corresponding empirical measure
	\begin{equation}
	p({\bf x}_{0:k}|{\bf y}_{1:k-1})\approx  \frac{1}{M}\sum_{i=1}^{M} \delta({\bf x}_{0:k} - {\bf x}^{(i)}_{0:k|k-1}).
	\end{equation}
	\textbf{Weight computation:}
	Importance weights are then computed according to an incoming observation ${\bf y}_k$:
	\begin{equation}
	\label{eq:weightsingle}
	w_{k|k}^{(i)} = \frac{p({\bf y}_k|{\bf x}^{(i)}_{k|k-1})}{\sum\limits_{l=1}^{M} p({\bf y}_k|{\bf x}^{(l)}_{k|k-1} )}  .
	\end{equation}
	
	\noindent\textbf{Resampling:}
	Lastly, one performs resampling by sampling with replacement the entire state trajectories ${\bf x}^{i}_{0:k|k-1}$, $i=1,\ldots,M$, according to the weighted 
	empirical measure
	\begin{equation}
	p({\bf x}_{0:k}|{\bf y}_{1:k})\approx\sum_{i=1}^{M} w_{k|k}^{(i)} \delta({\bf x}_{0:k}- {\bf x}^{(i)}_{0:k|k-1})
	\end{equation}
	instead of just the states at the current time (as is the case in filtering) \cite{Kitagawa1996}.  As a result one obtains $M$ equally weighted
	trajectories ${\bf x}_{0:k|k}^{(i)}$. \\ %\medskip \\
	
%\noindent
We remark that the resampling with replacement step can be thought of as applying an $M\times M$ transformation matrix ${\bf D}_{k|k}$ to the 
matrices ${\bf X}_{l|k-1}$ at each time instance $l=0,\ldots,k$, that is,
\begin{equation}  \label{eq:transform_l}
{\bf X}_{l|k} = {\bf X}_{l|k-1}{\bf D}_{k|k}.
\end{equation}
The entries of ${\bf D}_{k|k}$ are drawn at random from zero and one such that ${\bf D}_{k|k}^{\rm T} \mathbb{1} = \mathbb{1}$ and
\begin{equation}
\frac{1}{M} \mathbb{E}[{\bf D}_{k|k}\mathbb{1}] = {\bf w}_{k|k},
\end{equation}
where $\mathbb{1} \in \mathbb{R}^M$ denotes a vector of ones, 
\begin{equation}\label{eq:weights} 
	{\bf w}_{k|k}  = \left(w^{(1)}_{k|k},\ldots,w^{(M)}_{k|k} \right)^{\rm T} \in \mathbb{R}^{M\times 1}
\end{equation}
denotes the vector of importance weights (\ref{eq:weightsingle}), 
and $\mathbb{E}$ stands for taking the expectation value. If the $(i,j)$-th entry of ${\bf D}_{k|k}$ is one, then this implies
that the $i$-th ensemble member ${\bf x}_{l|k-1}^{(i)}$ at time $k$ gets replaced by the $j$-th member, that is ${\bf x}_{l|k}^{(i)} = 
{\bf x}_{l|k-1}^{(j)}$ for all $l=0,\ldots,k$.
	
	It is crucial to point out that a direct implementation of this classical particle smoother is known to degenerate as soon as $k$ is too large (even with resampling) \cite{saerkkae2013}. More specifically, given $M$ samples ${\bf x}_{0|0}^{(i)}$ from the initial distribution $p({\bf x}_0)$, collected in the ensemble
	matrix ${\bf X}_{0|0}$, one finds that the above algorithm leads to
	\begin{equation}
	{\bf X}_{0|k} = {\bf X}_{0|0} {\bf D}_{1|1} {\bf D}_{2|2} \cdots {\bf D}_{k|k}
	\end{equation}
	and ${\bf X}_{0|k}$ will contain many identical ensemble members, leading to a largely diminished effective sample size for the smoothing distribution 
	at initial time. Upon using the compact notation (\ref{eq:matrix_trajectories}) with $L = k$ one can rewrite the complete trajectory resampling step in 
	matrix form as
	\begin{equation} \label{eq:resample_T}
	{\bf X}_{0:k|k} = {\bf X}_{0:k|k-1}  {\bf D}_{k|k}
	\end{equation}
	from the columns of which one can read off the smoothed sample trajectories ${\bf x}_{0:k|k}^{(i)}$ if desired. 
	Again, many of these trajectories will contain identical state samples ${\bf x}_{l|k}^{(i)}$ as $l\ll k$.
	
	For this reason smoothing is often carried out for a relatively small fixed-lag $L$, that is, (\ref{eq:resample_T}) gets replaced by
	\begin{equation}  \label{eq:LETFanalysisupdate}
	{\bf X}_{k-L:k|k} = {\bf X}_{k-L:k|k-1}  {\bf D}_{k|k}
	\end{equation}
	with $L = 0$ leading back to filtering. Consequently we will focus on estimating the joint
	smoothing distributions $p({\bf x}_{k-L:k}|{\bf y}_{1:k})$ or the marginal distributions $p({\bf x}_l|{\bf y}_{1:k})$, $l = k-L,\ldots,k$, 
	in this paper.  Furthermore, recent transformation approaches to smoothing, such as the ESRS and the NETS, 
	can be put into the above framework (\ref{eq:LETFanalysisupdate}) with different choices of the matrix ${\bf D}_{k|k}$. Most importantly, 
	the entries of ${\bf D}_{k|k}$ are no longer drawn at random and can take real values.  We will also find that an extension of the 
	ETPF to smoothing requires an important modification to the above framework.

%%%%%%%%%%%%%%%%%%%%%%%%%%%%%%%%%%%%%%%%%%%%%%%%%%%%%%%%%%%%%%%%%%
	\section{Linear Ensemble Transform Smoothers}\label{sec:LETF}
%%%%%%%%%%%%%%%%%%%%%%%%%%%%%%%%%%%%%%%%%%%%%%%%%%%%%%%%%%%%%%%%%
We now introduce the concept of LETSs, using the framework developed in the previous section.  LETSs provide a generalisation of the associated family of linear ensemble transform filters 
(LETFs) \cite{sr:reichcotter15}.  An LETS is any smoother in which the ensemble of samples ${\bf X}_{k-L:k|k-1}$ of the prior $p({\bf x}_{k-L:k}|{\bf y}_{1:k-1})$ 
are moved towards an observation via a linear transformation in the form of (\ref{eq:LETFanalysisupdate}), where the transformation matrix must be of size $M \times M$.  
Note that $L=0$ corresponds to an LETF, and the transformation (\ref{eq:LETFanalysisupdate}) describes a standard filtering update in commonly employed filters such as the stochastic EnKF \cite{sr:evensen2006} and the ensemble square root filter (ESRF) \cite{sr:tippett03} \cite{sr:reichcotter15}. 

The remainder of this section is devoted to three specific examples of LETSs.  We first demonstrate that the ESRS and the NETS both fit into the framework provided by (\ref{eq:LETFanalysisupdate}). We then propose the smoother variant of the ETPF, showing that it is also an LETS, but that ${\bf D}_{k|k} $ in (\ref{eq:LETFanalysisupdate}) requires an important extension to ensure consistency.  

%%%%%%%%%%%%%%%%%%%%%%%%%%%%%%%%%%%%
	\subsection{Ensemble Square Root Smoother}
%%%%%%%%%%%%%%%%%%%%%%%%%%%%%%%%%%%
	 Smoother extensions of the class of ensemble square root filters \cite{sr:tippett03} have been investigated, for example \cite{Cosme2010, Fairbairn2009}.  Here we briefly summarize the extension of the ESRF in \cite{sr:cotterreich} in the context of LETSs.  The resulting smoother is a type of deterministic ensemble Kalman smoother, reformulated such that it can be seen as an LETS.
	 We assume for simplicity that the forward operator in (\ref{eq:t2}) is linear, that is, $h(x)={\bf H}x$. The ESRS is designed so that the empirical first two moments match the Kalman mean and covariance matrix under a linear model and observation setting. Then the coefficients of the linear transformation of the ESRS are given by
	\begin{equation} \label{eq:ESRS}
	{\bf D}^{\rm{ESRS}}_{k|k}  := \widehat{\bf w}_{k|k}\mathbb{1}^{\rm T} + \textbf{S}_{k|k}  ,
	\end{equation}
	with      
	\begin{equation}
	\widehat{\bf w}_{k|k} = \frac{1}{M-1}\textbf{S}_{k|k}^2 ({\bf H} \textbf{A}_{k|k-1})^{\rm T}\textbf{R}^{-1} ( {\bf y}_{k} - {\bf H}{{\bf m}}_{k|k-1} )
	\end{equation}
	and 
	\begin{equation}
	\textbf{S}_{k|k} = \left\{ {\bf I} + \frac{1}{M-1} ({\bf H} \textbf{A}_{k|k-1})^{\rm T} \textbf{R}^{-1} {\bf H} \textbf{A}_{k|k-1}\right\}^{-1/2},
	\end{equation}
	where the matrix of ensemble deviations and the ensemble mean are defined by
	\begin{equation}
	\textbf{A}_{k|k-1}= \Big({\bf x}^{(1)}_{k|k-1}-{{\bf m}}_{k|k-1},\dots,{\bf x}^{(M)}_{k|k-1}-{{\bf m}}_{k|k-1} \Big) 
	\end{equation}
	and
	\begin{equation}
	\label{eq:ensmean}
		{\bf m}_{k|k-1} =  {\frac{1}{M}\sum_{i=1}^{M} {\bf x}^{(i)}_{k|k-1}} ,
	\end{equation}
	respectively. The ESRS is obtained from (\ref{eq:LETFanalysisupdate}) with ${\bf D}_{k|k} = {\bf D}^{\rm ESRS}_{k|k}$. 
	Note that the transformation ${\bf D}^{\rm{ESRS}}_{k|k}$ is equal to the filter transformation of an ESRF at time $k$. In other words in an implementation of the ESRS one only needs to compute the filter transformations and multiply them with past states up to lag $L$.
	For more details on how to derive the transformation matrix (\ref{eq:ESRS}) see, for example, \cite{sr:reichcotter15}.
	
	The ESRS is easy to implement and applicable even under relatively small samples sizes.  However in general it is not consistent, that is, it does not generate samples of the true posterior distribution other than the case of linear Gaussian systems even in the limit of $M\to\infty$.

%%%%%%%%%%%%%%%%%%%%%%%%%%%%%%%%%%%%%%%%%	
	\subsection{Nonlinear Ensemble Transform Smoother}
%%%%%%%%%%%%%%%%%%%%%%%%%%%%%%%%%%%%%%%%%%

	Along the lines of the ESRS, the NETS \cite{Kirchgessneretal2017} is also constructed to match the first two moments correctly for finite $M$. The difference to the ESRS is that the NETS is matched to the empirical mean and covariance of the classical particle filter for fixed $M$.  This ansatz yields the following linear transformation
	\begin{equation}\label{NETSupdate}
	{\bf D}_{k|k}^{\rm{NETS}}= {\bf w}_{k|k}\mathbb{1}^{\rm T}+\sqrt{M}\left[{{\bf W}_{k|k}}-{\bf w}_{k|k}({\bf w}_{k|k})^{\rm T}\right]^{1/2}{\bf \Omega}	_{k|k}
	\end{equation}
	where ${\bf w}_{k|k}$ is the vector of importance weights (\ref{eq:weights}) at time $k$ and
	\begin{align}
	\label{eq:defnW}
{\bf W}_{k|k} = \mbox{diag}\,({\bf w}_{k|k})
	\end{align}
is a diagonal matrix with the importance weights as diagonal elements. The NETS is obtained from (\ref{eq:LETFanalysisupdate}) with ${\bf D}_{k|k} = {\bf D}^{\rm NETS}_{k|k}$. Again the transform matrix ${\bf D}^{\rm NETS}_{k|k}$ is equivalent to that of the corresponding NETF.
See \cite{Kirchgessneretal2017} for more details.
	
	The ${\bf \Omega}_{k|k}$ in (\ref{NETSupdate}) represents a suitably chosen orthogonal rotation matrix subject to ${\bf \Omega}_{k|k}\mathbb{1} = 
	\mathbb{1}$.  The choice of this matrix is discussed in further detail in Section \ref{sec:Adaptivespreadcorrectionandrotation}. Despite matching the first two moments of the particle smoother for a finite $M$ the filter does not converge towards the true posterior in the limit $M\to \infty$.

%%%%%%%%%%%%%%%%%%%%%%%%%%%%%%%%%%%%%%%%%%%%%%%%%%%%%%%%%	
	\subsection{Ensemble Transform Particle Smoother}\label{sec:ETPF}
%%%%%%%%%%%%%%%%%%%%%%%%%%%%%%%%%%%%%%%%%%%%%%%%%%%%%%%%%%
	
	The ETPF relies on a transformation which involves solving an optimal transport problem designed to maximize the correlation between $  {\bf X}_{k|k} $ and ${\bf X}_{k|k-1}$, whilst obtaining a set of samples distributed in accordance with the importance weights ${\bf w}_{k|k}$ \cite{sr:reichcotter15}.  The resulting transformation matrix, denoted by ${\bf D}^{\rm ETPF}_{k|k}$, is given by
	\begin{align}
				\label{eq:ETPFupdate}
				{\bf D}^{\rm ETPF}_{k|k} = \arg \min  \sum^M_{i,j=1} d_{ij}\,||{\bf x}^{(i)}_{k|k-1}-{\bf x}^{(j)}_{k|k-1}||^2
		 \end{align}
		 with $d_{ij}$ corresponding to the $(i,j)$-th element of the transformation matrix. 
	Following the formulation of the ESRS and the NETS, one could be tempted to define the ETPS extension of the ETPF by
setting ${\bf D}_{k|k} = {\bf D}^{\rm ETPF}_{k|k}$ in (\ref{eq:LETFanalysisupdate}). However such a naive approach does not take into account the
geometry of the underlying trajectory samples ${\bf x}_{k-L:k|k-1}^{(i)}$ and can lead to statistically inconsistent updates even in the limit $M\to \infty$
as demonstrated by the following simple experiment.
	
\begin{figure}
\begin{center}
	\includegraphics[width=8cm, height=6.857cm]{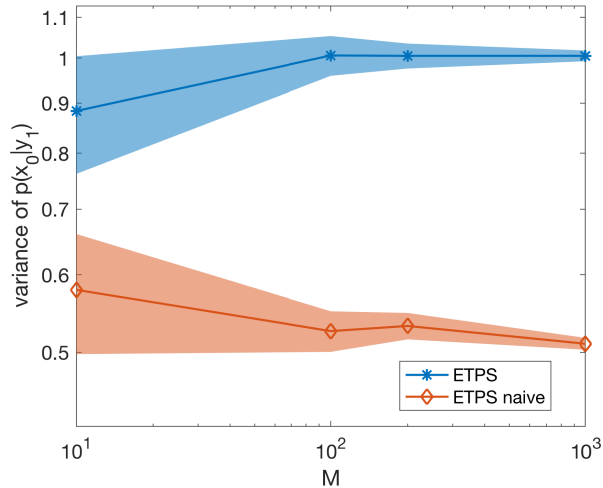} 
\end{center}
\caption{Estimated variance of the marginal $p({\bf x}_{0}|{\bf y}_{1})$ for $60$ independent runs for the ETPS update given in (\ref{eq:optimaltransportproblem}) and the naive extension of the ETPF (\ref{eq:ETPFupdate}).  The solid line shows the mean over the 60 experiments, and the shaded regions indicate the $95\%$ confidence intervals.}
\label{fig:Example1}
\end{figure}

\textbf{Example:} We consider the following one dimensional model with $f=0$ and ${\bf Q} = 1$ in (\ref{eq:t1}), i.e, 
	\begin{equation} \label{eq:ex1model}
	\begin{aligned}
	x_{k} = \epsilon_k, \ \epsilon_k \sim \mathcal{N}(0,1) 
	\end{aligned}
	\end{equation}
The observations are according to (\ref{eq:t2}) with ${\bf R}=1$ and $h(x) = x$, i.e,
\begin{equation} \label{eq:ex1data}
	\begin{aligned}
	y_{k} = h(x_{k}) + \nu_k,  \ \nu_k \sim \mathcal{N}(0,1) 
	\end{aligned}
	\end{equation}
	The purpose of this example is to demonstrate the inconsistency of the \underline{naive} approach and it is thus sufficient to only consider one filtering step (at time $t=1$). Let us assume that $y_1 = 0$. Note that $x_0\sim \mathcal{N}(0,1)$ and the prior distribution at time 1 is $p( x_{1})= \mathcal{N}(0,1)$ as it is independent of the state at $t=0$. The filtering distribution at $t=1$ is
	\begin{equation}\label{eq:ex1filteringdensity}
p(x_{1}|y_{1})=\mathcal{N}(0,0.5) 
	\end{equation}
whereas the marginal smoothing distribution at $t=0$ is
	\begin{equation}\label{eq:ex1marginalsmoothingdensity}
p(x_{0}|y_{1})=\mathcal{N}(0,1).
	\end{equation}
Note that (\ref{eq:ex1marginalsmoothingdensity}) is equal to the initial distribution $p(x_0)$ because $x_1$ and $x_0$ are independent of each other and thus the observation $y_1$ does not affect $x_0$. However, this is not what one obtains from the naive smoothing extension of the ETPF to this problem. Instead one finds that the variance of $p(x_0|y_1)$ is systematically reduced to approximately $0.5$ for $M\to \infty$. This is shown in Figure \ref{fig:Example1} which shows the variance of (\ref{eq:ex1marginalsmoothingdensity}) for $60$ individual independent runs for different ensemble sizes $M\in\{10,100,200,1000\}$ for the naive implementation of the ETPS, i.e., the update is done via the ETPF update given (\ref{eq:ETPFupdate}). For comparison, the variance obtain via the ETPS proposed below in $\ref{eq:optimaltransportproblem}$ is also displayed, which clearly converges to the correct variance as $M \rightarrow \infty$. 
%\end{example}

Instead the desired transformation matrix ${\bf D}_{k-L:k|k}^{\rm ETPS}$ should be obtained as the solution of the following linear transport problem:
\begin{equation}\label{eq:optimaltransportproblem}
	{\bf D}^{\rm ETPS}_{k-L:k|k}=\arg \min  \sum^M_{i,j=1} d_{ij}\,||{\bf x}^{(i)}_{k-L:k|k-1}-{\bf x}^{(j)}_{k-L:k|k-1}||^2
	\end{equation}
	with the minimization performed over the set of admissible $M\times M$ matrices ${\bf D}$ defined by
	\begin{subequations}\label{eq:constraintotimaltransportproblem}
	\begin{align}
 & d_{ij} \geq 0 \ \forall \ i,j \\
	 &{\bf D} \mathbb{1} = M{\bf w}_{k|k}  \\
	&{\bf D}^{\rm T} \mathbb{1} = \mathbb{1} .
\end{align}
\end{subequations}
Here $d_{ij}$ stands for the $(i,j)$-th entry of the matrix ${\bf D}$.
The ETPS is now obtained from (\ref{eq:LETFanalysisupdate}) with ${\bf D}_{k|k} = {\bf D}^{\rm ETPS}_{k-L:k|k}$. 

The transformation of the ETPS can be interpreted as a deterministic replacement of the resampling step of sequential importance resampling. The ETPS 
analysis step is consistent in the sample limit $M\to\infty$ as the transformation ${\bf D}^{\rm ETPS}_{k-L:k|k}$ converges to a transfer plan $\Psi$ inducing a deterministic coupling between the random variable associated with the prior trajectories ${\bf x}_{k-L:k|k-1}^{(i)}$ and the random variable associated with the posterior trajectories ${\bf x}_{k-L:k|k}^{(i)}$. This follows from a general approximation result of \cite{sr:mccann95}.

It is further important to note that the associated coupling is maximizing the correlation between the random variable associated with the prior and the one associated with the posterior.  Lastly, note that the ETPS transformation matrix in (\ref{eq:optimaltransportproblem}) is computed using both the current weights ${\bf w}_{k|k}$ and information regarding the ensemble of trajectories  ${\bf x}_{k-L:k|k-1}^{(i)}$ over the fixed-lag window. Hence the transform matrix 
${\bf D}^{\rm ETPS}_{k-L:k|k}$ is different from that of the corresponding ETPF at time $k$ which relies on the ensemble of states ${\bf x}_{k|k-1}^{(i)}$
only.

We demonstrate in Section \ref{sec:numerics}, that it is also preferable to replace ${\bf D}_{k|k}^{\rm NETS}$ in the NETS by a sample trajectory
based ${\bf D}_{k-L:k|k}^{\rm NETS}$ through an optimal choice of the rotation matrix $\mathbf{\Omega}_{k|k}$ in
(\ref{NETSupdate}). See also Section \ref{sec:Adaptivespreadcorrectionandrotation} below. Hence we close this section by noting that
\begin{equation}  \label{eq:LETS_proper}
	{\bf X}_{k-L:k|k} = {\bf X}_{k-L:k|k-1}  {\bf D}_{k-L:k|k}
\end{equation}
provides the appropriate general LETS framework with (\ref{eq:LETFanalysisupdate}) for ${\bf D}_{k-L:k|k} = {\bf D}_{k|k}$ as a special case.		
	
%%%%%%%%%%%%%%%%%%%%%%%%%%%%%%%%%%%%%%%%%%%%%%%%%%%%%
	\section{Numerical tools for moderate sample sizes}
%%%%%%%%%%%%%%%%%%%%%%%%%%%%%%%%%%%%%%%%%%%%%%%%%%%%%%

	Smoothing is naturally a higher dimensional problem than the associated filtering problem and thus it is important to address this additional 
	complexity in the context of moderate sample sizes. In other words, smoothers typically require additional techniques to make them feasible for practical 
	estimation problems. We will discuss three approaches in this section that can help to significantly improve the accuracy and computational 
	feasibility of the proposed family of smoothers especially in the small to moderate sample size setting.  We also discuss how the LETS framework is beneficial in such developments, since such techniques do not need to be developed for each algorithm individually, but can be applied to any smoother that can be seen as an LETS.  Further techniques for spatially extended
	systems will be discussed in Section \ref{sec:spatially_extended}.

%%%%%%%%%%%%%%%%%%%%%%%%%%%%%%%%%%%%%%%%%%%%%%
	\subsection{Temporal Localisation} \label{sec:TLocalization}
%%%%%%%%%%%%%%%%%%%%%%%%%%%%%%%%%%%%%%%%%%%%%%

Localisation has become a powerful tool for beating the curse of dimensionality for sequential filtering algorithms such as the ensemble Kalman filter.
These techniques have recently been extended to particle filters \cite{jdw:ChengReich,Poterjoy2015, Robert2017} 
and to smoothing algorithms. For example a weight localization procedure for use in Variational Particle Smoothers is proposed 
in \cite{Morzfeld2017}, where the weights are computed for each state component independently (or some block form of it).  \cite{Bocquet2016} 
investigate how a localization operator can be modified over the smoothing window using the underlying dynamics of the system.    

We note that the fixed-lag smoother can already been viewed as a temporal localisation of the global update step (\ref{eq:resample_T}). Still, 	
the optimal transport problem over an $(L+1)N_x$ dimensional space, as defined by (\ref{eq:optimaltransportproblem}) and (\ref{eq:constraintotimaltransportproblem}), becomes high-dimensional for time-lags $L \gg 1$ even if the dimension of state space, $N_x$, is moderate. Hence the ETPS might perform poorly under those circumstances even if the associated ETPF does well on a given problem. 
This degradation of performance can be counteracted by a further localization step 
in time. More precisely, we replace the single $(L+1)N_x$ dimensional transport problem by a sequence of $N_x$-dimensional  transport problems
	\begin{equation}\label{eq:optimaltransportproblem_l}
	{\bf D}^{\rm ETPS}_{l|k}=\arg \min  \sum^M_{i,j=1} d_{ij}\,||{\bf x}^{(i)}_{l|k-1}-{\bf x}^{(j)}_{l|k-1}||^2,
	\end{equation}
$l = k-L,\ldots,k$, with the minimization performed over the $M\times M$ matrices ${\bf D}$ subject to (\ref{eq:constraintotimaltransportproblem}). 
We then replace (\ref{eq:transform_l}) by
	\begin{equation} \label{eq:transform_lETPS}
	{\bf X}_{l|k} = {\bf X}_{l|k-1}{\bf D}^{\rm ETPS}_{l|k}
	\end{equation}
for $l = k-L,\ldots,k$. We note that the localized update (\ref{eq:transform_lETPS}) is now consistent with the ETPF filtering update since 
${\bf D}_{k|k}^{\rm ETPF} = {\bf D}_{k|k}^{\rm ETPS}$. However, while the marginal distributions $p({\bf x}_l|{\bf y}_{1:k})$ are correctly 
approximated by (\ref{eq:optimaltransportproblem_l}), the correlation between samples ${\bf x}_{l|k}^{(i)}$ at different time-levels 
$l,l' \in \{k-L,\ldots,k\}$ can be misrepresented, which can lead to problems when applied to a model of the form (\ref{eq:t1b})  of order 
$n$, $n>1$, where the joint distribution of state variables over a time interval $[k-n, k-1]$ matters in the prediction step. 
We explore this issue in Section \ref{sec:numerics} in more detail. 

Having to solve $L+1$ optimal transport problems (\ref{eq:optimaltransportproblem_l}) can become computationally demanding. 
Our numerical experiments for the Lorenz 63 system from Section \ref{sec:numerics} 
suggest that this cost can be reduced by using 
\begin{equation}
\label{eq:compcheaperloc}
{\bf D}_{l|k}^{\rm ETPS} = {\bf D}_{k|k}^{\rm ETPS} \quad \text{for all} \quad l=k-L,\ldots,k, 
\end{equation}
provided the sample paths are strongly correlated in time, as is the case for many deterministic systems (i.e. with ${\bf Q} = 0$).  This issue is investgated further in Section \ref{sec:numerics}.  It should be mentioned that the update formulation of LETSs in (\ref{eq:LETFanalysisupdate}) involves modifying the ensemble at each time $l$ in the smoothing window.  When applied to deterministic systems, the resulting ensembles will not necessarily satisfy the underlying system of equations. This issue is well known in the area of ensemble Kalman filtering/smoothing, and methods to ensure system balances are maintained is an active area of research.  One could also adopt an alternative approach whereby the goal is to estimate the smoothing update at the beginning of the window only, similar to the iEnKS \cite{Bocquet2014} and others.

%%%%%%%%%%%%%%%%%%%%%%%%%%%%%%%%%%%%%%%%%%%%
	\subsection{Hybrid Smoothers} \label{sec:Hybrid}
%%%%%%%%%%%%%%%%%%%%%%%%%%%%%%%%%%%%%%%%%%%
	Ensemble smoothers with an underlying Gaussian assumption are often favoured over consistent smoothers, as they are typically more robust in high dimensional applications. 
	Yet depending on the problem, a Gaussian approximation can lead to a very crude estimate of the true posterior. One way to improve accuracy while maintaining some level of robustness is to combine different smoothers. This strategy has proven to be rather successful in the filtering context \cite{sr:frei13,Potthastetall2019}.
	Our hybrid formulations are based on the idea of a split likelihood \cite{sr:CRR15}:
	\begin{equation}\label{splitlikelihood}
	p({\bf y}_{k}|{\bf x}_{k}) = p({\bf y}_{k}|{\bf x}_{k})^\alpha \,p({\bf y}_{k}|{\bf x}_{k})^{(1-\alpha)}
	\end{equation}
	where $\alpha \in [0,1]$ is a tuning parameter that determines the contribution of each smoother.  The benefit of the LETS framework is that it allows one to develop a hybrid smoother using any pair of smoothers that can be written in the form of (\ref{eq:LETFanalysisupdate}).  Specifically, this is done by a two-step transformation to form a mixture that benefits from desirable properties of each individual smoother.  For example, a hybrid of the ETPS and ESRS would be developed as
	\begin{equation}
	\label{eq:eghybrid}
	{\bf X}_{k-L:k|k} = {\bf X}_{k-L:k|k-1}{\bf D}^{\rm{ETPS}}_{k-L:k|k}(\alpha) {\bf D}^{\rm{ESRS}}_{k|k} (1-\alpha).
	\end{equation}
	Here ${\bf D}^{\rm{ETPS}}_{k-L:k|k}(\alpha)$ denotes the
	transformation matrix corresponding to the tempered likelihood $p({\bf y}_{k}|{\bf x}_{k})^\alpha$ and ${\bf D}^{\rm{ESRS}}_{k|k} (1-\alpha)$ to
	$p({\bf y}_{k}|{\bf x}_{k})^{(1-\alpha)}$, respectively. Note that $p({\bf y}_{k}|{\bf x}_{k})^\alpha$, $\alpha \in (0,1]$, for given likelihood 
	(\ref{eq:likelihood}) corresponds to replacing the measurement error covariance matrix ${\bf R}$ by ${\bf R}/\alpha$ in the ESRS and ETPS formulations.
	The hybrid example in (\ref{eq:eghybrid}) is particularly valuable as two LETSs with different strengths are combined yielding a filter superior to the individual LETSs (for optimal $M$ and $\alpha$) \cite{sr:CRR15}. It has also been found numerically that lower RMSEs are achieved when the 
	ETPF is used first \cite{sr:CRR15}.  The choice of $\alpha$ determines the influence of each smoother, with either $0$ or $1$ indicating purely one of the smoothers.  The optimal $\alpha$ can be tuned using an appropriate performance metric, for example, the RMSE.

%%%%%%%%%%%%%%%%%%%%%%%%%%%%%%%%%%%%%%%%%%%%%%%%%%%%%%%
	\subsection{Adaptive spread correction and rotation}\label{sec:Adaptivespreadcorrectionandrotation}
%%%%%%%%%%%%%%%%%%%%%%%%%%%%%%%%%%%%%%%%%%%%%%%%%%%%%%%
	In practice, LETSs are often employed in a setting far from the asymptotic limit ($M\rightarrow\infty$) which often results in an underestimation of the spread. This underestimation is typically counteracted by the use of artificial inflation of the sample spread (sometimes also referred to as particle rejuvenation). These procedures are usually tied to a tuning factor that needs to be adjusted.
	
	An alternative approach is to correct the sample spread adaptively by matching it to the second moment of the classical particle smoother in order to improve the overall accuracy of LETSs.  
	This ansatz is justified by the fact that the classical particle smoother converges weakly to the true posterior and thus can be a good estimate of the true measure (if $M$ is not too small).  We describe this approach below.  Any smoother belonging to the subfamily of LETSs (\ref{eq:LETS_proper}) 
	with transformation matrices ${\bf D}_{k-L:k|k} \in \mathcal{D}_1$,
	\begin{equation}
	\mathcal{D}_1=\{{{\bf D}} 
	\in\mathbb{R}^{M\times M}|\, {{\bf D}}^{\rm T}\mathbb{1} = \mathbb{1},~ {{\bf D}}\mathbb{1} = M{{\bf w}_{k|k}}\, \},
	\end{equation}
	match the first moment of the particle filter over a single assimilation step.  Here the vector of weights ${\bf w}_{k|k}$ is defined by
	(\ref{eq:weights}).
	
	For example, the ETPS as well as the smoother associated with the simple transformation ${{\bf D}}_{k-L:k|k} = {\bf w}_{k|k}{\mathbb{1}}^{\rm T}$
	are in $\mathcal{D}_1$. Provided ${\bf D}_{k-L:k|k} \in\mathcal{D}_1$, it follows that this leads to equivalence of the empirical means, that is,
	\begin{equation}
	{\bf m}_{k-L:k|k}= {\bf X}_{k-L:k|k-1}{\bf w}_{k|k} = \frac{1}{M} {\bf X}_{k-L:k|k-1}{\bf D}_{k-L:k|k}\mathbb{1} = 
	\frac{1}{M} {\bf X}_{k-L:k|k}\mathbb{1}. 
	\end{equation}
	However smoothers (\ref{eq:LETS_proper}) with transformation ${\bf D}_{k-L:k|k} \in \mathcal{D}_1$ 
	do not necessarily match the second moment, that is,
	\begin{equation}\label{eq:conimportancesampling}
	{\bf P}_{k-L:k|k} = \sum_{i=1}^M w^{(i)}_{k|k-1} ({\bf x}^{(i)}_{k-L:k|k-1}-{\bf m}_{k-L:k|k})
	{({\bf x}^{(i)}_{k-L:k|k}-{\bf m}_{k-L:k|k})^{\rm T}}
	\end{equation} 
	is not equal to 
	\begin{equation}\label{eq:contransformation}
	\widehat{\bf P}_{k-L:k|k} = \frac{1}{M} \sum_{i=1}^M ({\bf x}^{(i)}_{k-L:k|k}-{{{\bf m}}_{k-L:k|k}) ({\bf x}^{(i)}_{k-L:k|k}-{{\bf m}}_{k-L:k|k})^{\rm T}}.
	\end{equation} 
	
	For (\ref{eq:conimportancesampling}) and (\ref{eq:contransformation}) to be equal the following equation has to be satisfied
	\begin{equation}
	({\bf D}_{k-L:k|k}  - {\bf w}_{k|k}\mathbb{1}^{\rm T})({\bf D}_{k-L:k|k}  - {\bf w}_{k|k}\mathbb{1}^{\rm T})^{\rm T} = M( 
	{\bf W}_{k|k} - {\bf w}_{k|k}{\bf w}_{k|k}^{\rm T}),
	\end{equation}
	where ${\bf W}_{k|k} $ is defined in (\ref{eq:defnW}).  \\
\textbf{Remark} It is important to note that for small ensemble sizes the spread given by (\ref{eq:conimportancesampling}) associated with importance sampling is expected to be underestimated as well. Yet ultimately the ETPS is underestimating the spread even stronger and in the context of filtering it has been confirmed numerically that an adjustment to the particle filter significantly improves the results even for a small number of particles \cite{acevedo2016}. Thus an adjustment of the ETPS to the particle smoother spread can be a beneficial modification even when far from the ensemble limit.

	The transformation can be corrected via an additive term ${\bf {\Delta}}_{k|k}$ that vanishes in the particle limit $M\to \infty$ similar to the adjustment in the filtering case \cite{acevedo2016}. That is, for any ${{\bf D}_{k-L:k|k}}\in\mathcal{D}_1$ the associated second order correct transformation is computed as
	\begin{equation}
	\widehat{{\bf D}}_{k-L:k|k} = {{\bf D}}_{k-L:k|k} + {\bf \Delta}_{k-L:k|k}
	\end{equation}
	where  ${\bf \Delta}_{k-L:k|k} \in \mathbb{R}^{M\times M}$ is chosen such that ${\bf \Delta}_{k-L:k|k}^{\rm T}\mathbb{1} = {\bf 0}$, 
	${\bf \Delta}_{k-L:k|k} \mathbb{1} = {\bf 0}$, and is a solution of the algebraic Riccati equation
	\begin{equation}
	\begin{split}
	\label{riccati}
	&M ({\bf W}_{k|k} - {\bf w}_{k|k}{\bf w}_{k|k}^{\rm T} ) 
	-({{\bf D}}_{k-L:k|k} - {\bf w}_{k|k}\mathbb{1}^{\rm T})
	({{\bf D}}_{k-L:k|k} - {\bf w}_{k|k}\mathbb{1}^{\rm T})^{\rm T}  \\ &
	\qquad  = 
	({{\bf D}}_{k-L:k|k}-{\bf w}_{k|k}\mathbb{1}^{\rm T}) {\bf \Delta}_{k-L:k|k}^{\rm T} + {\bf \Delta}_{k-L:k|k} 
	({{\bf D}}_{k-L:k|k} - {\bf w}_{k|k}\mathbb{1}^{\rm T})^{\rm T}
	+ {\bf \Delta}_{k-L:k|k}{\bf \Delta}_{k-L:k|k}^{\rm T}.
	\end{split}
	\end{equation}
	Such Riccati equations are in general non-trivial to solve, however for the purposes of adjusting the transformation it is often sufficient to approximate the solution by numerically solving the corresponding dynamic Riccati equation \cite{acevedo2016}.  The above developments demonstrate the usefulness of the LETS framework, as one only needs the form of the transformation matrix (as well as the particle smoother weights) in order to construct a second order correction.
	Further, it is important to note that  
	\begin{equation}\label{eq:correctedversionwithrotation}
	{\bf D}_{k-L:k|k}^{\bf \Omega} = {\bf D}_{k-L:k|k} + {\bf \Delta}_{k-L:k|k}{\bf \Omega}_{k-L:k|k},
	\end{equation}
	where ${\bf \Omega}_{k-L:k|k}$ is an $M\times M$ orthogonal matrix with the property that ${\bf \Omega}_{k-L:k|k}\mathbb{1} = \mathbb{1}$ is also a second order accurate transformation since $ {\bf \Delta}_{k-L:k|k}{\bf \Omega}_{k-L:k|k}$ is also a solution to (\ref{riccati}). This fact is particularly relevant as the second order corrected version of ${{\bf D}}_{k-L:k|k} = {\bf w}_{k|k}{\mathbb{1}}^{\rm T}$ leads to the NETS (\ref{NETSupdate}) and the correction term takes the form
	\begin{equation} \label{eq:delta}
	{\bf \Delta}_{k-L:k|k} = \sqrt{M}\left[{{\bf W}_{k|k}}-{\bf w}_{k|k}({\bf w}_{k|k})^{\rm T}\right]^{1/2}.
	\end{equation}
	In \cite{acevedo2016} it was shown how the performance of the NETF ($L=0$) can be improved
	significantly by a suitably chosen ${\bf \Omega}_{k-L:k|k}$. This approach can be extended to the smoothing problem 
	whereby ${\bf \Omega}_{k-L:k|k}$ is chosen so that ${\bf D}_{k-L:k|k}^{\bf \Omega}$ is the minimizer of 
	 \begin{equation}\label{eq:NETSoptimaltransport}
	  \mathcal{V}({\bf D}^{\bf \Omega}) = \sum^M_{i,j=1} d^{{\bf \Omega}}_{ij}\,||{\bf x}^{(i)}_{k-L:k|k-1}-{\bf x}^{(j)}_{k-L:k|k-1}||^2 , \qquad {\bf D}^{\bf \Omega} = {\bf w}_{k|k}{\mathbb{1}}^{\rm T}
	  + {\bf \Delta}_{k-L:k|k}{\bf \Omega}
	 \end{equation}	
over the set of suitable rotation matrices ${\bf \Omega}$ for ${\bf \Delta}_{k-L:k|k}$ given by (\ref{eq:delta}).	 

The desired rotation matrix can be determined via the singular value decomposition of an $M\times M$ matrix, that is,
\begin{align}\label{eq:ensembledeviations}
{\bf \Delta}_{k-L:k|k}^{\rm T} {\bf A}_{k-L:k|k-1}^{\rm T} {\bf A}_{k-L:k|k-1}  
={\bf U}{\bf \Lambda} {\bf V}^{\rm T}.
\end{align}
with the $(L+1)N_x \times M$ matrix of ensemble deviations defined by
\begin{equation}
{\bf A}_{k-L:k|k-1} := {\bf X}_{k-L:k|k-1} - \frac{1}{M} {\bf X}_{k-L:k|k-1}\mathbb{1} \mathbb{1}^{\rm T}.
\end{equation}
Then the optimal rotation in (\ref{eq:correctedversionwithrotation}) for ${\bf D}_{k-L:k|k} = {\bf w}_{k|k}{\mathbb{1}}^{\rm T}$ is given by
\begin{equation}
{\bf \Omega}_{k-L:k|k} ={\bf U} {\bf V}^{\rm T}.
\end{equation}
The reader is referred to \cite{acevedo2016} for further details.

In Section \ref{sec:numerics}, we demonstrate how an optimal rotation of this form can significantly improve the accuracy of the NETS with random rotation introduced in \cite{Kirchgessneretal2017}.

Two fixed lag smoothers belonging to the family of LETSs, namely the ETPS and the NETS with optimal rotation, are summarized 
in form of pseudocodes below.  More specifically, the ETPS given in Algorithm \ref{alg:ETPFsmoother} combines the ETPS of Section 
\ref{sec:ETPF} with the additional feature of a second order correction as described above. 
In order to compute a solution to the optimal transport problem given in step (\ref{eq:optimaltransportLag}) of Algorithm \ref{alg:ETPFsmoother}, one can employ available solvers such as the Earth Mover's Distances algorithm given in \cite{sr:Pele-iccv2009}. The associated computational complexity is of order $\mathcal{O}(M^3\ln(M))$. For the theoretical aspects associated with optimal transport problems please see volumes such as \cite{sr:Villani,sr:Villani2} and references therein. In the context of filtering the optimal transport problem has been treated in \cite{sr:reichcotter15}.
As it is computationally expensive to solve the optimal transport problem it is possible to utilize the so called Sinkhorn approximation (see \cite{sr:cuturi13}) of the actual solution ${\bf D}^{\rm{ETPS}}_{k-L:k|k}$ with computational complexity of order $\mathcal{O}(M^2 C(\lambda))$. The Sinkhorn approximation involves adding a regualrisation term to the underlying transport problem, i.e,
\begin{equation}
\label{eq:sinkhorn}
{\bf D}^{\rm{ETPS}}_{k-L:k|k}(\lambda) = \arg \min \sum_{i,j=1}^M \left\{ d_{ij} || {\bf x}_{k-L:k|k-1}^{(i)} - {\bf x}_{k-L:k|k-1}^{(j)} ||^2 + \frac{1}{\lambda} d_{ij} \log \left(\frac{d_{ij}}{d_{ij}^0} \right) \right\}
\end{equation}
where $\lambda > 0$ is a regularisation parameter that must be tuned and $d_{ij}^0$ are the entries of ${\bf D}^0_{k-L:k|k} = {\bf w}_{k|k} \mathbb{1}^T.$  Note that this form of approximation does not ensure that the smoother is consistent in the ensemble limit anymore. Further it is important to point out that the an update with a transformation computed via the Sinkhorn approximation underestimates the spread even more than with the original ETPS transformation. Thus it is prudent to pair the Sinkhorn approximation with the the adaptive spread correction.
Although the computational complexity of solving (\ref{eq:optimaltransportLag}) is not significantly affected by the smoother lag $L$ as the transformation matrix has dimension $M\times M$ regardless of the value of $L\ge 1$, the computational burden of solving the optimal transport problem is significant. Thus an accuracy complexity trade-off in form of the Sinkhorn approximation is justifiable. For a more detailed discussion in the context filtering we refer to \cite{acevedo2016}.

The NETS has already been introduced in \cite{Kirchgessneretal2017} and the novelty in the disclosed pseudocode (see Algorithm \ref{alg:NETS}) is the use of an optimal rotation matrix $ {\bf \Omega}_{k-L:k|k}$ in (\ref{eq:correctedversionwithrotation}) for ${\bf D}_{k-L:k|k} = {\bf w}_{k|k} 
\mathbb{1}^{\rm T}$ and ${\bf \Delta}_{k-L:k|k}$ given by (\ref{eq:delta}).  

\begin{algorithm}[h]
	\caption{ETPS with fixed-lag}
	\label{alg:ETPFsmoother}
	\begin{algorithmic}[1]
		\State{Set: \begin{itemize}
				\item ensemble size $M$
				\item smoother lag $Lag$
			\end{itemize}
		}
		\State{Initialize ${\bf x}^{(i)}_0\sim p({\bf x}_0)\quad i=1,\dots,M$}
		\For{$k=1:T$}
		\State\label{line3}{Generate ${\bf x}^{(i)}_{k|k-1}\sim p({\bf x}_k|{\bf x}^{(i)}_{k-1|k-1})$,
		update the ensemble matrix ${\bf X}_{k-L:k-1|k-1}$ by augmentation to ${\bf X}_{k-L:k|k-1}$,  
		and determine normalized importance weights according to (\ref{eq:weightsingle}).}
		\State{$L = \min(k,Lag)$}
		\State\label{Recomputationtransformationmatrix}{Compute transformation matrix via
			\begin{equation*}\label{eq:optimaltransportLag}
			{\bf D}^{\rm{ETPS}}_{k-L:k|k}=\arg \min  \sum^M_{i,j=1} d_{ij}\,||{\bf x}^{(i)}_{k-L:k|k-1}-{\bf x}^{(j)}_{k-L:k|k-1}||^2
			\end{equation*}
			subject to constraints (\ref{eq:constraintotimaltransportproblem}).}
		\State{For second-order accurate ETPS: calculate ${\bf\Delta}_{k-L:k|k}$ determined by (\ref{riccati}) with ${\bf D}_{k-L:k|k} = 
		{\bf D}^{\rm{ETPS}}_{k-L:k|k}$, otherwise: let ${\bf\Delta}_{k-L:k|k} = {\bf 0}$}  \label{secondorderlineAlg}  
		\State{Define  
			\begin{equation}
			 {\bf  D}^{\rm{SETPS}}_{k-L:k|k}  ={\bf D}^{\rm{ETPS}}_{k-L:k|k}+ {\bf \Delta}_{k-L:k|k}
			\end{equation}}
		\State{Update 
			\begin{equation}
			{\bf X}_{k-L:k|k}  = {\bf X}_{k-L:k|k-1}   {\bf  D}^{\rm{SETPS}}_{k-L:k|k} 
			\end{equation}}
			\EndFor
		 \end{algorithmic} 
		\end{algorithm}

\begin{algorithm} 
\caption{optimal NETS with fixed-lag}
\label{alg:NETS}
\begin{algorithmic}[1]
\State{Set: \begin{itemize}
\item ensemble size $M$
\item smoother lag $Lag$
\end{itemize}
}
\State{Initialize ${\bf x}^{(i)}_0\sim p({\bf x}_0)\quad i=1,\dots,M$ }
\For{$k=1:T$}
	\State\label{line4}{Generate ${\bf x}^{(i)}_{k|k-1}\sim p({\bf x}_k|{\bf x}^{(i)}_{k-1|k-1})$, update the ensemble matrix 
	 ${\bf X}_{k-L:k-1|k-1}$ by augmentation to ${\bf X}_{k-L:k|k-1}$,  and determine normalized importance weights 
	${\bf w}_{k|k}$ according to (\ref{eq:weightsingle}).}
	\State{$L = \min(k,Lag)$}
	\State{Compute 
	\begin{equation}
	{\bf \Delta}_{k-L:k|k} =\sqrt{M} ({\bf W}_{k|k} - {\bf w}_{k|k}{\bf w}_{k|k}^{\rm T})^{1/2}
	\end{equation} }
		\State{
		Compute singular value decomposition of (\ref{eq:ensembledeviations}) and define ${\bf \Omega}_{k-L:k|k} = {\bf U} {\bf V}^{\rm T}$.
		}
\State{Define  
			\begin{equation}
			 {\bf  D}^{\rm{NETS}}_{k-L:k|k}  = {\bf w}_{k|k} \mathbb{1}^{\rm T} + {\bf \Delta}_{k-L:k|k} {\bf \Omega}_{k-L:k|k}
			\end{equation}}
		\State{Update
			\begin{equation}
			{\bf X}_{k-L:k|k}  = {\bf X}_{k-L:k|k-1}   {\bf  D}^{\rm{NETS}}_{k-L:k|k} 
			\end{equation}}
			
\EndFor
\end{algorithmic} 
\end{algorithm}

%%%%%%%%%%%%%%%%%%%%%%%%%%%%%%%
%
\section{Numerical examples} \label{sec:numerics}
%
%%%%%%%%%%%%%%%%%%%%%%%%%%%%%%%
In order to complement the above derivations we conducted three sets of numerical experiments each targeting a different investigational purpose.

%%%%%%%%%%%%%%%%%%%
\subsection{Classical: Lorenz 63}\label{se:L63}
%%%%%%%%%%%%%%%%%%%%

\begin{figure}
	\begin{center}
		\includegraphics[width=16cm, height=13.714cm]{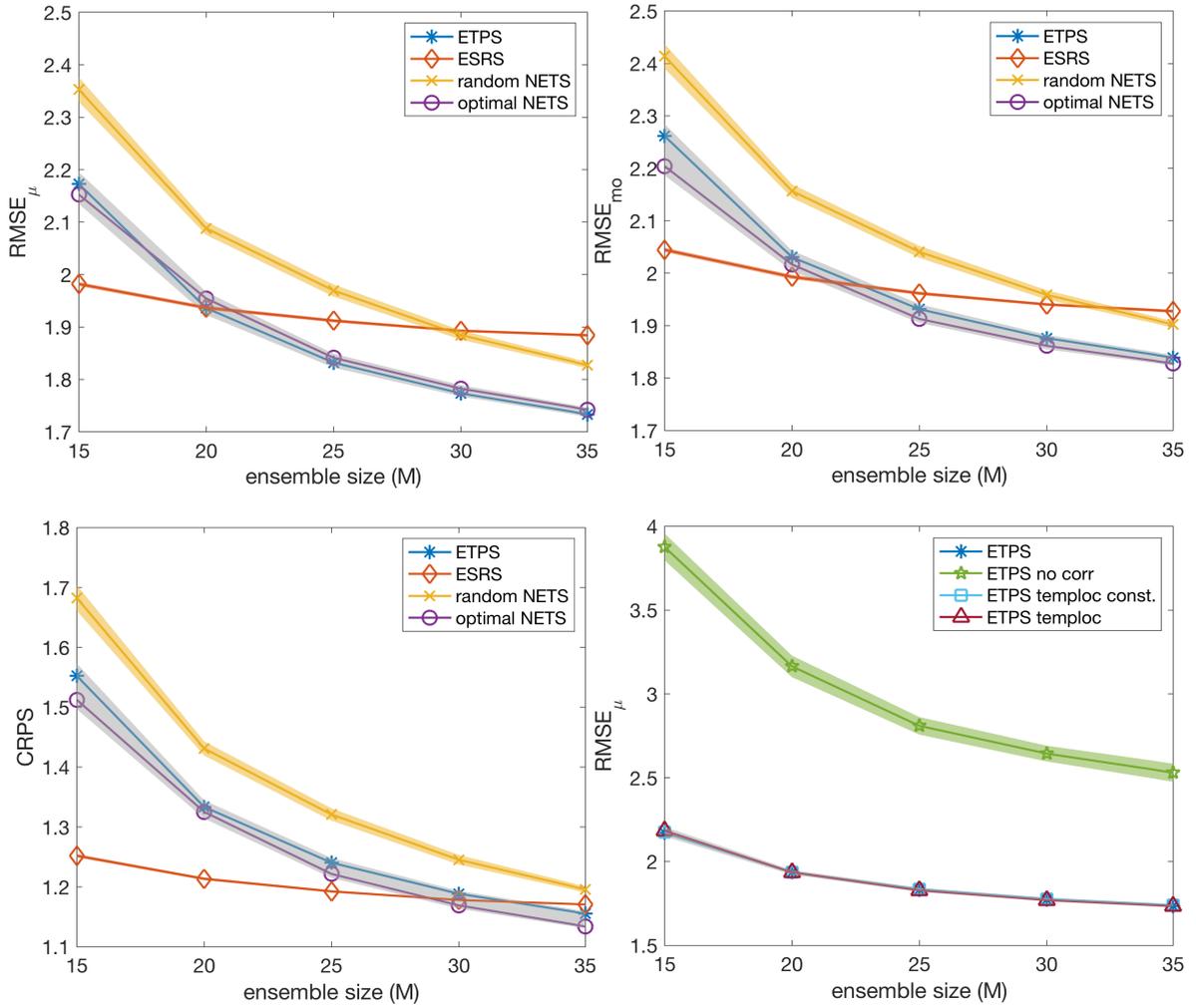}
		
		\caption{Summary statistics vs ensemble size for various smoothers at optimal fixed lag $L=6$  for the Lorenz 63 experiment.  Time averaged RMSE computed on ensemble mean (top left) and ensemble mode (top right) as well as CRPS (bottom left) are shown, in addition to various modifications of the ETPS (bottom right).  \textcolor{black}{Case temploc refers to the localisation scheme in (\ref{eq:optimaltransportproblem_l}), while temploc const. corresponds to (\ref{eq:compcheaperloc}).  For reference, a bootstrap particle smoother with $M=2000$ and $L=6$ gave the following summary statistics: $\text{RMSE}_{\mu} = 1.2, \text{RMSE}_{mo} = 1.29, CRPS = 0.69$.  Case no corr refers to the ETPS without second order correction discussed in Sec \ref{sec:Adaptivespreadcorrectionandrotation}.}   Optimal NETS refers to the NETS with the proposed rotation (Algorithm \ref{alg:NETS}) whilst random rotation refers to the rotation scheme in \cite{Kirchgessneretal2017}.  Shaded regions indicate the $95 \%$ confidence intervals.  The gray shaded region indicates the region that bounds all $95 \%$ confidence intervals of the ETPS and optimal NETS, same for the bottom right subplot except for the ETPS and the two localisation schemes.}
		\label{fig:L63_RMSE_optimallag}
	\end{center}
\end{figure}

The Lorenz 63 system \cite{sr:lorenz63} is a deterministic ordinary differential equation in $\mathbb{R}^3$ with ${\bf x}_k = (x_k, y_k, z_k)^T$ given by
\begin{equation} 
\frac{d{\bf x}_k}{dt}= \psi({\bf x}_k) := \begin{bmatrix}      
10(y_k - x_k)\\
x_k(28 - z_k) - y_k\\
x_k y_k - 8/3 z_k
\end{bmatrix}.
\end{equation}
Forward Euler time-stepping with step-size $\Delta t$ leads to a system of the 
form (\ref{eq:t1}) with $f({\bf x}_k) = {\bf x}_k + \Delta t \psi ({\bf x}_k)$ and ${\bf Q} = {\bf 0}$.
It is a classical toy example in the context of data assimilation as it is chaotic and solutions diverge rapidly without any filtering updates. Here we mainly want to compare some of the LETSs equipped with the various implementation variants proposed in this paper. In detail a standard implementation of the ESRS, NETS with different rotations (random \cite{sr:toedter15} and optimal as proposed in Algorithm \ref{alg:NETS} based on the ideas in \cite{acevedo2016}) and the ETPS (with and without second order correction see right panel of Figure \ref{fig:L63_RMSE_optimallag}). Further a hybrid smoother combining the ESRS with the ETPS is applied to the system.
For the computations we consider the case where only the first component is observed, at intervals of $\Delta t_{obs} = 0.12$ with observation error variance $R = 8$.  The discretised system is evaluated using $\Delta t = 0.01$.  Additionally, the initial conditions are uncertain, with an ensemble generated from $N({\bf x}_0, 0.5{\bf I})$.  Due to the deterministic model dynamics, particle rejuvenation is carried out as follows:

\begin{equation}
\label{eq:rejuv}
{\bf x}_{k|k}^{(j)} \rightarrow {\bf x}_{k|k}^{(j)} + \beta \left( \frac{1}{M-1}\sum_{i=1}^{M} ({\bf x}_{k|k-1}^{(i)} - {\bf m}_{k|k-1}) ({\bf x}_{k|k-1}^{(i)} -{\bf m}_{k|k-1})^{\rm T}\right)^{1/2} {\bf \xi}^{(j)}  
\end{equation}
where $\beta$ is a tuning parameter taken to be 0.2 in this experiment and ${\bf \xi}_j$ are i.i.d Gaussian variables with mean 
${\bf \mu} = {\bf 0}$ and covariance matrix ${\bf I}$.  This choice of parameters has been checked to be robust to small variations and is a reasonable choice compared to values used in other studies.

A sequence of $K = 10,000$ observations is generated and assimilation carried out using 50 identical experiments for each algorithm, due to the randomness associated with rejuvenation. Smoother performance for a given fixed lag $L$ is assessed using the time averaged RMSE:
\begin{equation}
\label{eq:RMSE}
\text{RMSE}(L) = \frac{1}{K}\sum_{k=1}^{K}\sqrt{\frac{1}{N_x}[{\bf \hat{x}}_{k-L|k} - {\bf x}^{\text{ref}}_{k-L}][{\bf \hat{x}}_{k-L|k} - {\bf x}^{\text{ref}}_{k-L}]^{\rm T} } %, \quad k-L \leq l < k, 
\end{equation}
where $K$ is the length of assimilation period (i.e., RMSE is evaluated only at time points where observations are available); ${\bf x}^{\text{ref}}_{k-L}$ refers to the true state at time $k-L$ and ${\bf \hat{x}}_{k-L|k}$ is a summary statistic of the smoothed ensemble at time $k-L$ given data up to time $k$.  We consider the traditional RMSE computed on the ensemble mean, denoted $\text{RMSE}_\mu$ with ${\bf \hat{x}}_{k-L|k} = {\bf m}_{k-L|k}$.  Since the densities in the nonlinear examples in this section may be skewed or display multi-modality, we also assess performance by a RMSE as per (\ref{eq:RMSE}), but computed on the sample mode rather than the ensemble mean ${\bf m}_{k-L|k}$, denoted $\text{RMSE}_{mo}$.  The sample mode is estimated from the maximum of a kernel density estimate fitted to the marginal ensembles at each time.  Finally, we examine the time and space-averaged Continuous Ranked Probability Score (CRPS) \cite{Hersbach2000} for each fixed lag $L$ to examine reliability and resolution of the posterior ensembles.

The main interest here is the performance of various smoothers for small to moderate ensemble sizes.  Figure \ref{fig:L63_RMSE_optimallag} shows the time averaged RMSEs and CRPS for the fixed lag $L = 6$, where RMSE is calculated at observation time points only.  The optimal fixed lag (i.e. with lowest $\text{RMSE}_\mu$) was 5 or 6 across all ensemble sizes $M = 15$ to 35.  In this experiment, 1 lag unit is equal to 1 assimilation interval.  For the ETPS, larger values of lamba were shown to provide improved RMSE values (results are shown for $\lambda = 40$, and $\lambda = 100$ for the ETPS and ETPS with no second order correction respectively). Note that the choice of $\lambda$ also affects run time, i.e., large $\lambda$ values increase run time but also improve the approximation of the true optimal transport problem. In a real world application, an ideal lambda would need to be determined based on a desired balance between runtime and accuracy. As the number of ensemble members is not greater than $M=40$ we do not choose $\lambda>40$. The only exception is the ETPS run without spread correction in the right panel of Figure \ref{fig:L63_RMSE_optimallag} where the regularisation is chosen to be $\lambda=100$ to indicate that even for larger $\lambda$, one is not able to obtain good RMSEs without the spread correction.

\begin{figure}[h]
	\begin{center}
		\includegraphics[width=8cm, height=6.857cm]{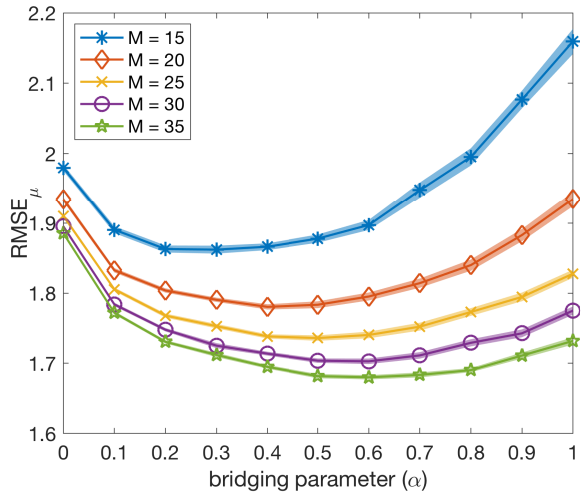} 
		\caption{Time averaged $RMSE_{\mu}$ vs bridging parameter $\alpha$ for the hybrid second-order accurate ETPS-ESRS for optimal fixed lag $L = 6$ for the Lorenz 63 experiment. The ETPS corresponds to $\alpha = 1$ and ESRS corresponds to $\alpha = 0$. } 
		\label{fig:L63_RMSE_hybrid}
	\end{center}
\end{figure}

Firstly, this experiment demonstrates the potential of the rotation techniques discussed in Section \ref{sec:Adaptivespreadcorrectionandrotation}. The NETS with optimal rotation provides improved results compared to the random rotation proposed in \cite{Kirchgessneretal2017}, to the extent that performance is similar to the 2nd order corrected ETPS.  Consistent with filtering results presented in \cite{acevedo2016}, the ESRS also provides the best smoother performance for low ensemble sizes, although the ETPS and NETS provide lower RMSEs for larger ensemble sizes both in terms of the ensemble mean and mode.  Such improvements can be attributed to the fact that the ETPS and NETS are better equipped to handle nonlinear relationships between observed and unobserved variables, due to the use of importance sampling rather than the linear correlation in the Kalman update.  However, there is a delicate balance between bias and variance here.  Larger ensemble sizes are required for the NETS and ETPS to become comparable to the ESRS in terms of CRPS.  This can be attributed to the slight underdispersiveness in the ETPS and NETS ensembles in this experiment.  Similar to the filtering case, the ETPS without second order correction gives the worst performance of all the smoothers.  
We also consider a hybrid smoother formulation as described in Section \ref{sec:Hybrid}, where a smoother update is first computed using the second order accurate ETPS which is then used as the prior for an update using the ESRS. Figure \ref{fig:L63_RMSE_hybrid} shows $\text{RMSE}_\mu$ values for the ETPS-ESRS hybrid smoother at optimal lag $(L = 6)$ when the ETPS update is computed prior to the ESRS update.  Improved $\text{RMSE}_\mu$ values compared to both the standard ETPS $(\alpha = 1)$ and ESRS $(\alpha = 0)$ are achieved in the hybrid framework. 

\subsubsection{Optimal lag}  \label{sec:optimallag}

\begin{figure}
	\begin{center}
		\includegraphics[width=8cm, height=6.857cm]{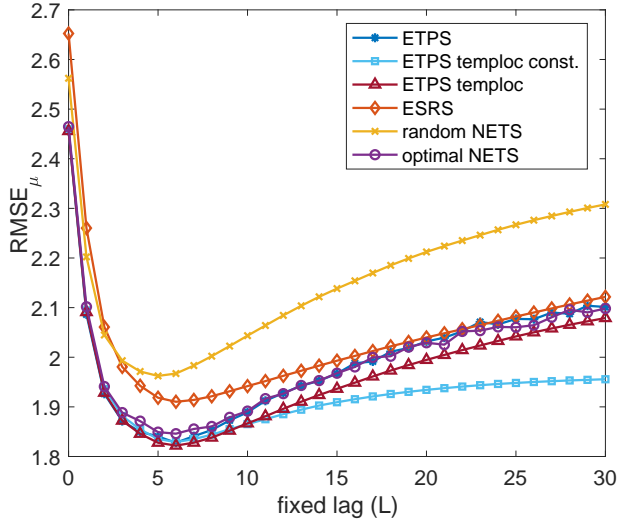} 
		\caption{Time averaged $\text{RMSE}_{\mu}$ vs fixed lag for various smoothers and localisation schemes for $M = 25$ for the Lorenz 63 experiment.  \textcolor{black}{Case temploc refers to the localisation scheme in (\ref{eq:optimaltransportproblem_l}), while temploc const. corresponds to (\ref{eq:compcheaperloc}).}   Optimal NETS refers to the NETS with the proposed rotation (Algorithm \ref{alg:NETS}) whilst random rotation refers to the rotation scheme in \cite{Kirchgessneretal2017}.  Note that 1 lag unit = 1 assimilation interval. }
		\label{fig:L63_RMSE_lag}
	\end{center}
\end{figure}

\begin{figure}[h]
	\begin{center}
		\includegraphics[width=16cm, height=6.857cm]{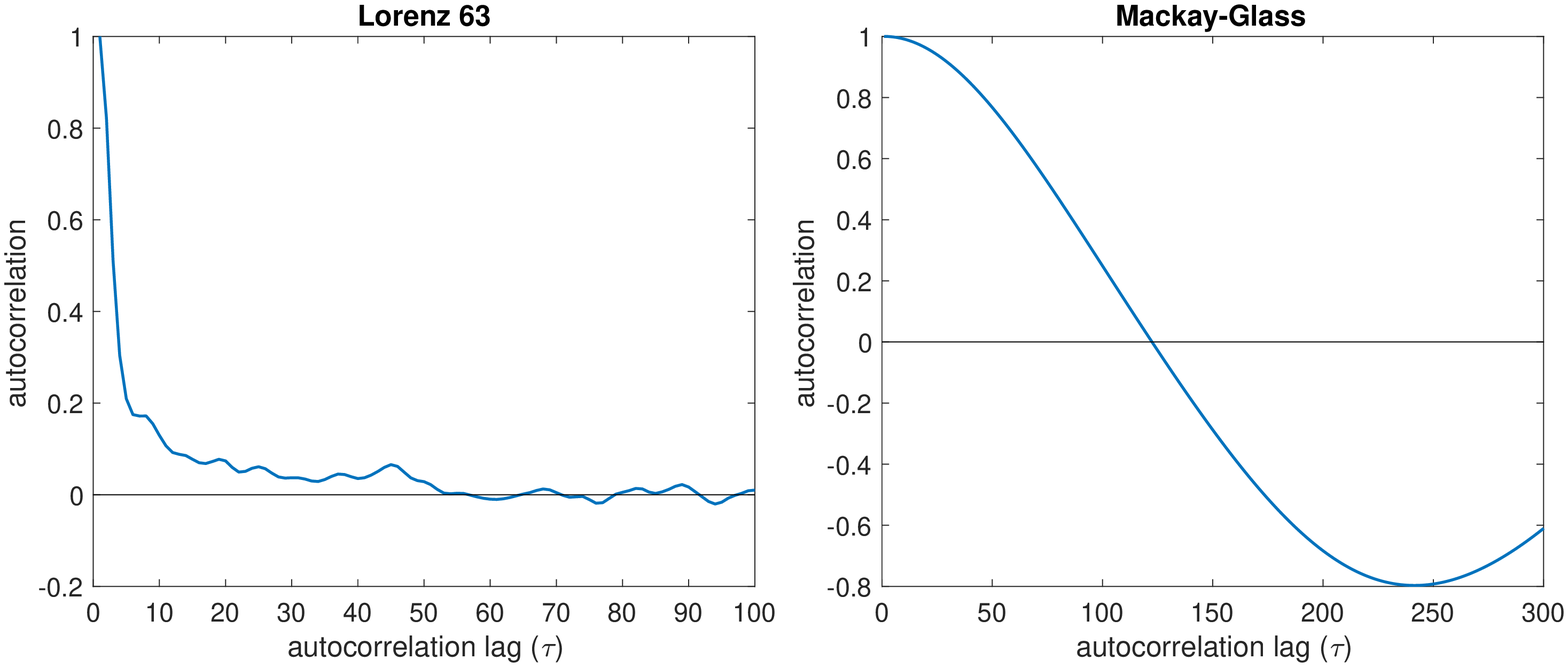} 
		\caption{Autocorrelation function for the x-component (i.e. observed component) of the Lorenz 63 system (left) and the Mackey-Glass system (right).  Note that for the Lorenz 63, 1 lag unit = 1 assimilation interval; for the Mackey-Glass system, 1 lag unit = 1 model time step.}
		\label{fig:ACF_MG_L63}
	\end{center}
\end{figure}

Due to sampling errors, it is often found that the optimal smoother performance is obtained for finite lags.  This phenomenon is often observed for the Lorenz systems when the ensemble size is small (see \cite{Kirchgessneretal2017} for an example with the Lorenz 96 system \cite{sr:lorenz96}) and also Figure \ref{fig:L63_RMSE_lag}  where the $\text{RMSE}_\mu$ increases for $L > 6$.  This may also be related to the autocorrelation properties of the underlying system; if the autocorrelation is small or negligible after a finite time, then it is possible that smoothing beyond this time lag will have minimal benefits.  These autocorrelation characteristics may be affected by the underlying chaotic properties of the system and how perturbations propagate in time.  Note the correspondence between the optimal smoother lag for the Lorenz 63 $(L = 5$ to $6)$ and the autocorrelation lag $(\tau)$ beyond which the autocorrelation function (ACF) of the $x$ variable  is small $( < 0.2)$ (see Figure \ref{fig:ACF_MG_L63}).  It is worth noting that the ACF provides only a measure of linear correlation and therefore cannot be used alone to determine the optimal smoother lag.  The potential of the localisation schemes can also be seen in Figure \ref{fig:L63_RMSE_lag} (right). The extreme localisation scheme of (\ref{eq:compcheaperloc} (see temp loc const.) provides significantly improved $\text{RMSE}_\mu$  at longer time lags.  We note that in practice, one would not continue to smooth past the optimal lag time and that this example is merely for demonstration purposes.

%%%%%%%%%%%%%%%%%%%%%%%%%%%%%%%%%
\subsection{Long Memory: Mackey-Glass model} \label{sec:MGM}
%%%%%%%%%%%%%%%%%%%%%%%%%%%%%%%%%%

	An important application of smoothing is in non-Markovian systems.  These systems typically have strong temporal correlations in the state variables, meaning that fixed lag smoothing over a long window has the potential to improve state estimation compared to filtering.  To investigate this further, we examine the Mackey-Glass model (originally proposed by \cite{Mackey1977}, which is a single-variable non-linear delay differential equation used to model various physiological processes such as blood cell production.  It is particularly useful for such systems because the evolution of the state is known to depend on the state at previous times, not just the current time.  The system has infinite degrees of freedom since it is a delay differential equation, however its strange attractor is of finite dimension.  Here, we consider the following parameterization of the model which produces periodic and chaotic solutions
\begin{equation}
\frac{dx}{dt}=\phi \frac{x(t-\nu)}{1+x(t-\nu)^\kappa}-\gamma x(t) \quad \text {with } \nu, \gamma,\phi,\kappa>0
\end{equation}
with $\phi = 0.2, \gamma = 0.1, \kappa = 10$ and $\nu = 17$. We use a 4th order Runge-Kutta scheme with $\Delta t = 0.1$ to discretise the system, resulting in a finite difference system in the form of (\ref{eq:t1b}) with $n = \frac{\nu}{\Delta t} = 170$.  Observations of the state variable are available at intervals of $\Delta t_{obs} = 8$, with observation error variance $R = 0.05$.  Here we consider the case where observations are temporally sparse since the system consists of only a single state variable.  The initial condition samples $x_{0:n}^{(i)}$ are generated by perturbing $x_0$ with random noise from $\mathcal{N}(0,0.1)$ and then propagating it through the model until time index $n$.  As with the Lorenz 63 system, particle rejuvenation using (\ref{eq:rejuv}) with $\beta = 0.005$\footnote{This choice is robust with respect to small changes.} is carried out since the model dynamics are deterministic.  A sequence of $K = 12,000$ observations is generated and assimilation carried out using 50 identical experiments for each smoothing algorithm.  The ETPS was solved using the Sinkhorn approximation (\ref{eq:sinkhorn}) with $\lambda = 40$ (see discussion in L63 section). In this experiment we consider only the NETS with optimal rotation proposed in Section \ref{sec:Adaptivespreadcorrectionandrotation}.    

\begin{figure}
	\begin{center}
		\includegraphics[width=16cm, height=13.714cm]{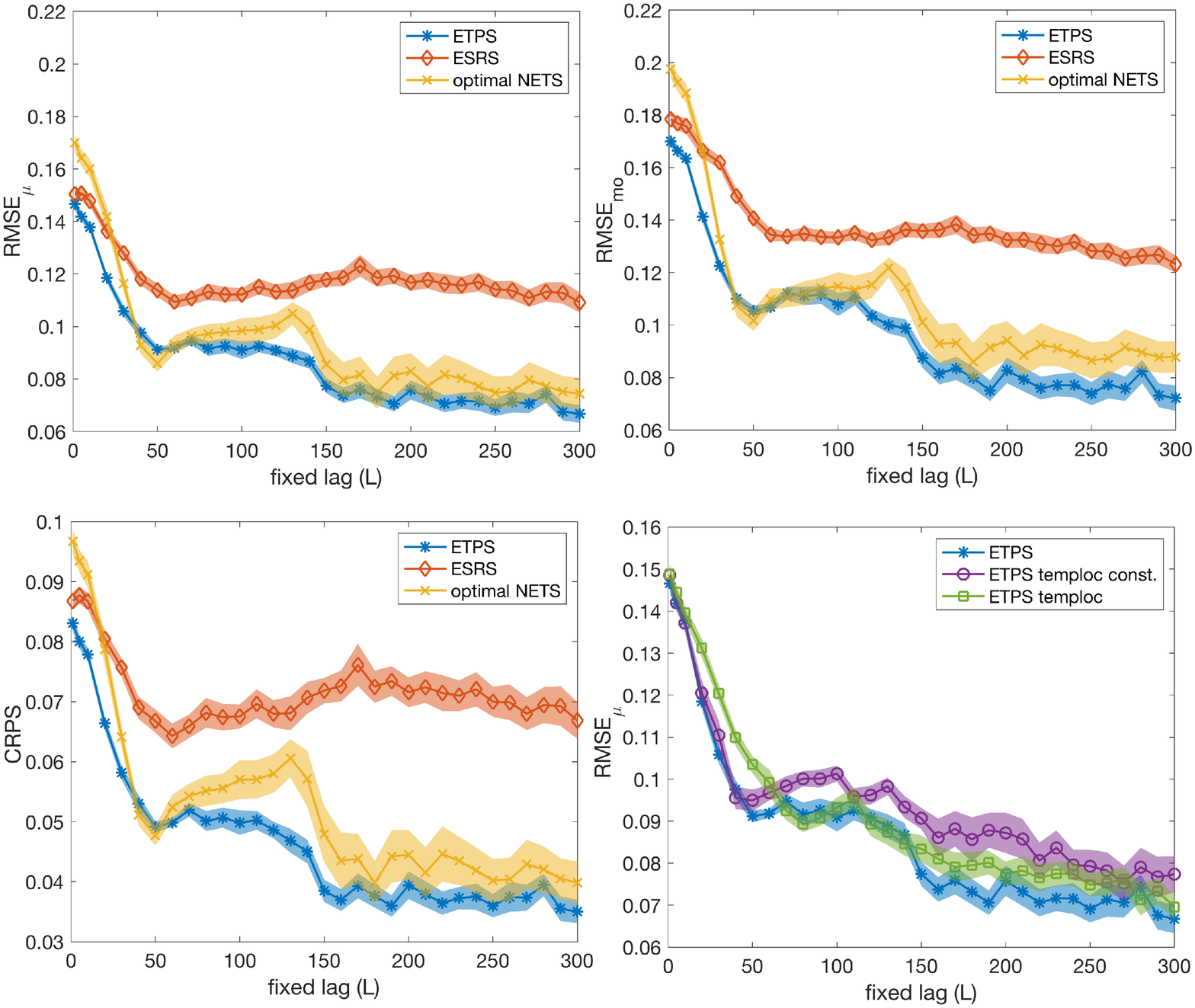}
		\caption{Summary statistics vs fixed lag for various smoothers for $M = 50$ for the Mackey-Glass experiment.  Time averaged RMSE computed on ensemble mean (top left) and ensemble mode (top right) as well as CRPS (bottom left) are shown, in addition to various localization schemes for the ETPS (bottom right).  Case temploc refers to the scheme in (\ref{eq:optimaltransportproblem_l}), while temploc const. corresponds to (\ref{eq:compcheaperloc}). Shaded regions indicate the $95 \%$ confidence intervals.  Note that 1 lag unit = 1 model time step.   Optimal NETS refers to the NETS with the proposed rotation (Algorithm \ref{alg:NETS}).   }
		\label{fig:RMSE_lag_MG}
	\end{center}
\end{figure}

Figure \ref{fig:RMSE_lag_MG} shows the time averaged RMSEs and CRPS of the state variable against fixed lag $L$ for various smoothers for $M = 50$.  In this experiment, 1 lag unit is equal to 1 model time step, since the system is higher order Markov.  RMSE is also evaluated at every time step rather than only when observations are available, so that we can evaluate performance at the hidden states also.  Once again, the ETPS algorithm provides a significant improvement over the ESRS, and a slight improvement over the NETS with the optimal determinstic rotation scheme, in terms of all performance metrics.  There is also an overall trend in reduction of all metrics with increasing fixed lag size, unlike in the Lorenz 63 system.  This can be attributed to the presence of strong temporal correlations in the system even at long lag times (see Figure \ref{fig:ACF_MG_L63}).  Figure \ref{fig:RMSE_lag_MG} (right) shows that the temporal localisation schemes discussed in Section \ref{sec:TLocalization} leads to only a minor loss of accuracy, with the computationally cheaper scheme of (\ref{eq:compcheaperloc}) still providing competitive results compared to the pathwise ETPS and localisation scheme of  (\ref{eq:optimaltransportproblem_l}).  Finally, note that RMSE does not decrease monotonically with lag because the presence of a higher order memory term in the model equation means that smoother solutions for lag $L$ are not solely a deterministic transformation of smoother solutions at lag $L-1$.   

%%%%%%%%%%%%%%%%%%%%%%%%%%%%%%%%%%%%%%%%%%%%%%%
\section{Spatially extended systems} \label{sec:spatially_extended}
%%%%%%%%%%%%%%%%%%%%%%%%%%%%%%%%%%%%%%%%%%%%%%%%
In many applications the state space of the considered systems can be very high dimensional which can be challenging if the ensemble size is considerably smaller, i.e., $N_x\gg M$. Localisation techniques have been successfully employed in these settings \cite{sr:hunt07,sr:houtekamer01,tNerger2014} to address this impacts of small ensemble sizes with an extended state space. In Section \ref{sec:TLocalization} we already highlighted one form of temporal localisation for the ETPS which has been numerically investigated in Sections \ref{sec:MGM} and \ref{se:L63}. Here we will further explore localisation and look into options that allow to localise jointly in space and time for the various smoothers mentioned in this manuscript. Moreover an atmospheric toy model often used to investigate the feasibility of newly developed methods for extended state spaces is considered.

\subsection{Temporal and Spatial Localisation} \label{sec:SLocalization}
%%%%%%%%%%%%%%%%%%%%%%%%%%%%%%%%%%%%%%%%%%%%%%%%%%%
In the following we restrict our attention to problems where the state vector is defined on a spatial domain that is a set of discrete points in the Euclidean space $\mathbb{R}^n$, denoted by $\{{\bf p}_j\}_{j = 1:N_x}$.  That is, the $s$-th component of the vector ${\bf x}_k$, denoted ${\bf x}_k [s]$, gives the value of some variable at the position ${\bf p}_s \in \mathbb{R}^n$.   For example, this may be a scalar PDE in one space dimension (i.e., ${\bf p} \in \mathbb{R}$) or on a two dimensional grid. Similarly, the observation vector is defined on a spatial domain in the same Euclidean space, but potentially on a different set of discrete points $\{{\bf q}_j\}_{j = 1:N_y}$.  We also assume that the observation locations are fixed in time and that observation errors are independent in space and time (i.e. ${\bf R}$ is diagonal). Furthermore, we use the notation ${\bf z}[i]$ to indicate the $i$-the component of the vector ${\bf z}$ while ${\bf z}(r)$ indicates that the vector ${\bf z}$ is a function of the variable $r$.\\
The key idea behind localisation is that spurious temporal correlations \cite{tNerger2014} stemming from sampling errors due to small ensemble sizes relative to the dimension of the state space are artificially reduced \cite{sr:evensen2006}. 
Typically localisation is implemented by either adjusting the empirical covariance matrix (which is referred to as ${\bf B}$-localisation in the context of NWP \cite{sr:houtekamer01, Hamill2001}) or by localising in observation space, known as ${\bf R}$-Localisation. Direct adjustment of the covariance matrix is not extendable to particle filters thus we will focus on localisation in observation space. \\
The general motivation behind ${\bf R}$-Localisation is that most physical systems exhibit a spatial and temporal decay in correlation. It is thus pertinent to determine whether there is a decay in correlation, and if so, how variables are jointly correlated across time and space. This is clearly highly dependent on the characteristics of the dynamical system, so it is crucial to introduce a distance measure $d_l({\bf p}_s,{\bf q}_u)$ between the points ${\bf p}_s$ and ${\bf q}_u$ depending on time $l$ tailored to the underlying physical model. 

Distance in space as well as spatial correlations between state and observation components are often relatively easily to determine for many systems of interest. However temporal correlations are much harder to identify as they depend on the dynamical flow of the system, making localisation for smoothers more difficult. Yet there have been several suggestions in recent literature on how to approach this problem. For example, in \cite{Bocquet2016} the temporal localisation function is governed by a Liouville equation for the iterative ensemble Kalman smoother \cite{Bocquet2014}, while a weight localisation technique \cite{Poterjoy2015} based on the distance of the current state to the observation is presented for a variational particle smoother in \cite{Morzfeld2017}. A scheme that allows to adaptively incorporate how the covariance is affected by the flow is discussed in \cite{BishopHodyss2011} in form of ${\bf B}$-localisation.

Here we suggest two choices of distance measures for ${\bf R}$-Localisation that could be used in the context of smoothing. As the optimal-lag is often relatively short, one ansatz is to assume that the spatio-temporal dependence does not change significantly within the smoothing interval. With this idea in mind, one option is to keep the metric stationary in time i.e.,
\begin{equation}\label{eq:dfix}
d_l({\bf p}_s,{\bf q}_u) = m({\bf p}_s,{\bf q}_u) \quad \forall \enskip l = k-L, \cdots, k
\end{equation}  
where $m(\cdot, \cdot)$ refers to an appropriate distance metric, for example, the Euclidean distance, which also incorporates the underlying spatial boundary conditions of the system.
This ansatz is similar to the localisation function chosen in \cite{Morzfeld2017}. An alternative approach is to take the spatio-temporal dynamics into account in the localisation function. One way of doing so is to consider the autocorrelation structure of the underlying system.  Given a time series of model simulations $\{\mathcal{X}_t\}_{t=1:T}$, recall that the sample autocorrelation matrix  at autocorrelation lag $\tau$ is defined by 
	\begin{equation}\label{eq:autocorrcoef}
	{\bf \Gamma}_\tau = {\bf V}^{-1} {\bf \Sigma_\tau} {\bf V}^{-1} , \quad {\bf \Gamma}_\tau \in \mathbb{R}^{N_x \times N_x}
	\end{equation}
	where ${\bf \Sigma_\tau} $ is the sample autocovariance matrix given by  
	\begin{equation}
	{\bf \Sigma_\tau} = \frac{1}{T-1} \sum_{t = \tau+1}^T (\mathcal{X}_t - \hat{\mu}) (\mathcal{X}_{t-\tau} - \hat{\mu})'
	\end{equation}
	with $\hat{\mu} = \frac{1}{T} \sum_{t=1}^T \mathcal{X}_t$, and ${\bf V}$ is a diagonal matrix with the $i$-th entry on the diagonal given by $\sqrt{{\bf \Sigma}_0[i,i]}$.

Our proposed localisation scheme is to take the distance between any given observation point ${\bf q}_u$ and the spatial grid point ${\bf p}_{s^*}$ at time $k$, which is the grid point that is maximally correlated in magnitude with ${\bf p}_s$ at time $l$, i.e.

\begin{equation}\label{eq:dauto}
d_l({\bf p}_s,{\bf q}_u) = m({\bf p}_{s^*}, {\bf q}_u)
\end{equation}

	\begin{equation}
	s^*= \arg \max_i  \left| {\bf \Gamma}_{\tau} [s,i]  \right|.
	\end{equation}
	The smoothing lag variable is directly related to the autocorrelation lag, i.e. $l \equiv k- \tau$, where $\tau \in \{0, 1, \cdots, L\}$.  Furthermore, notice that (\ref{eq:dauto}) is independent of $k$. 

The next step of ${\bf R}$-Localisation is to adjust ${\bf R}^{\rm{loc}}_k(g_{sl})^{-1}$ according to the influence of the observation ${\bf y}_{k}$ on the individual components of the state ${\bf x}_{k-L:k}[g_{sl}]$ for all $s = 1, \cdots, N_x$ and $l = k-L, \cdots, k$, where $g_{sl}$ is the index of the vector ${\bf x}_{k-L:k}$ corresponding to the $s$-th position coordinate at the $l$-th time in the smoothing window $k-L$ to $k$. More specifically
\begin{equation}
{\bf R}^{\rm{loc}}_k(g_{sl})^{-1} = C_k(g_{sl}) {\bf R} ^{-1},
\end{equation}
where $C_k(g_{sl})$ is an $N_y\times N_y$ diagonal matrix with zeroes on the off-diagonal elements (since measurement errors are assumed independent) with the $u$-th diagonal element equal to 
\begin{equation}
\rho\Big(\frac{d_l({\bf p}_s,{\bf q}_u)}{r_{loc}}\Big)
\end{equation}
for all $u = 1, \cdots N_y$. $\rho$ is a localisation function (e.g., Gaspari-Cohn \cite{sr:gaspari99}) and $r_{loc}$ is the localisation radius which must be specified a priori. For the NETS and the ETPS the localised inverse error covariance ${\bf R}_k^{\rm{loc}}(g_{sl})^{-1}$ enters through the unnormalized weights 

\begin{equation}\label{likelihoodloc}
\tilde{w}_{k|k}^{(i)}(g_{sl})\propto\exp\Big(-\frac{1}{2}(h({\bf x}_{k})-{\bf y}_{k})^{\rm T}({\bf R}_k^{\rm{loc}}(g_{sl})^{-1}(h({\bf x}_{k})-{\bf y}_{k})\Big) 
\end{equation}
for all $s = 1, \cdots, N_x$ and $l = k-L, \cdots, k$. This form of localisation is an extension of localisation in the ETPF \cite{sr:cotterreich,jdw:ChengReich}. There have been several other suggestions on how to localise particle filters \cite{Poterjoy2015} as well as variants of the particle smoother \cite{Morzfeld2017}.
Note that localisation does not only help to mitigate particle degeneracy, which generally occurs as soon as the smoothing interval is too large, but it can also help to reduce run time. This is because the update for R-localisation is done individually for each component (or for smaller blocks) and can thus be parallelised. Furthermore, note that the ensemble update via ${\bf \tilde{D}}_{l|k}(g_{sl}) $ is now done individually for each component $g_{sl}$ of the state vector ${\bf x}_{k-L:k|k}$ which can mean reduced computational complexity for certain smoothers.  For example, in the case of the ETPS, the aforementioned localisation scheme is combined with the localisation scheme in Section \ref{sec:TLocalization}, which results in solving the following one dimensional optimal transport problem

\begin{equation}
	{\bf \tilde{D}}^{\rm ETPS}_{l|k}(g_{sl}) =\arg \min  \sum^M_{i,j=1} \tilde{d}_{ij}\,||{\bf x}^{(i)}_{k-L:k|k-1}(g_{sl}) - {\bf x}^{(j)}_{k-L:k|k-1} (g_{sl}) ||^2,
	\end{equation}
	subject to 
	\begin{subequations}
		\begin{align}
		& \tilde{d}_{ij} \geq 0 \ \forall \ i,j \\
		&{\bf D} \mathbb{1} = M{\bf \tilde{w}}_{k|k}(g_{sl})  \\
		&{\bf D}^{\rm T} \mathbb{1} = \mathbb{1} .
		\end{align}
\end{subequations}  
This can be solved via a sorting algorithm which considerably decreases the run time \cite{sr:cotterreich}.  Similarly, for the ESRS, such a localisation scheme means that the smoother can be parallelised since the update is carried out for each ${\bf x}_{k-L:k|k}[g_{sl}]$ independently.  Consequently localisation is often deliberately employed when dealing with high-dimensional state spaces and is thus a crucial feature to make smoothing feasible.   

%%%%%%%%%%%%%%%%%%%%%%%%%
\subsection{Lorenz 96}
\begin{figure}
	\begin{center}
			\includegraphics[width=8cm, height=6.857cm]{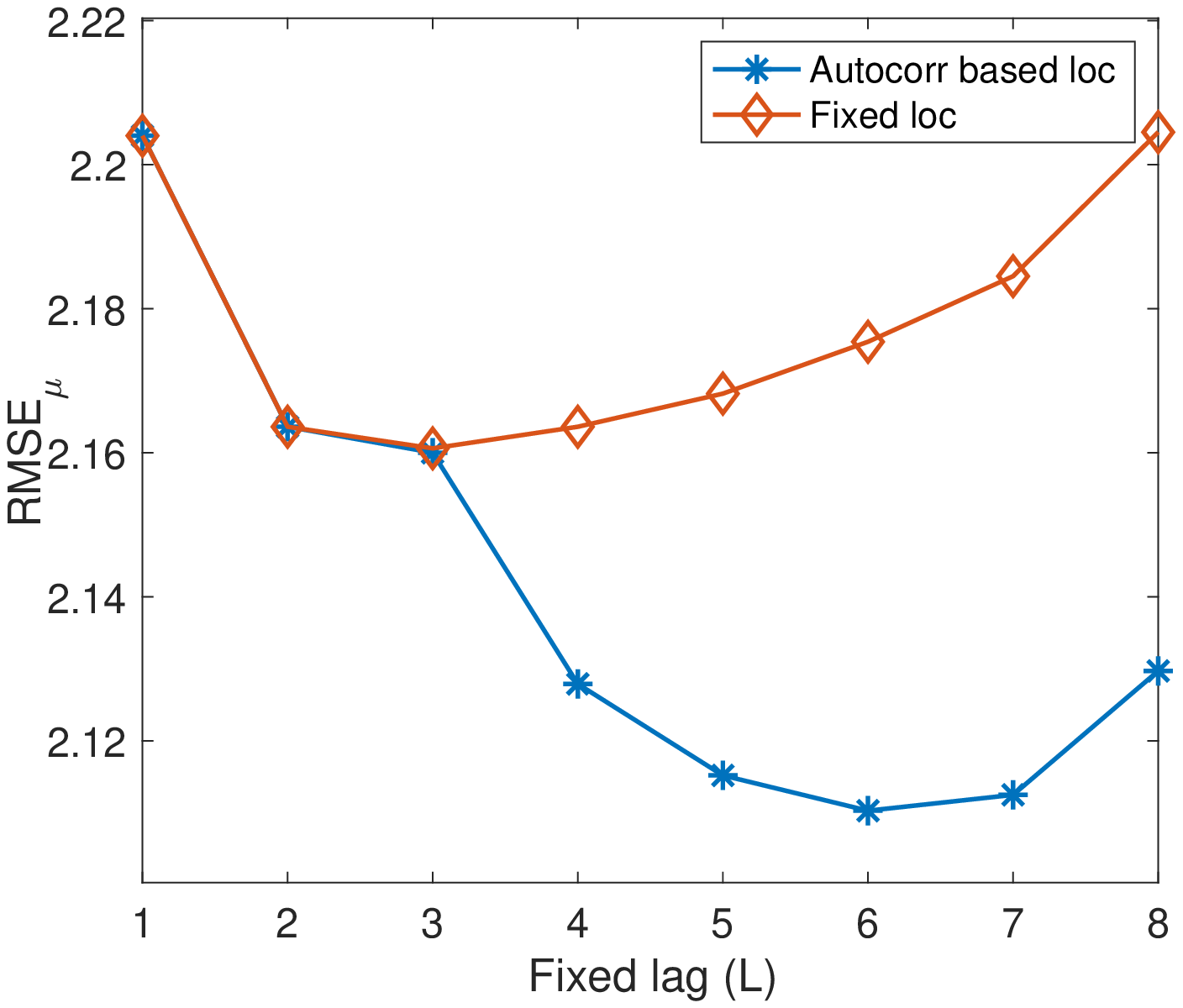} 
		\includegraphics[width=8cm, height=6.857cm]{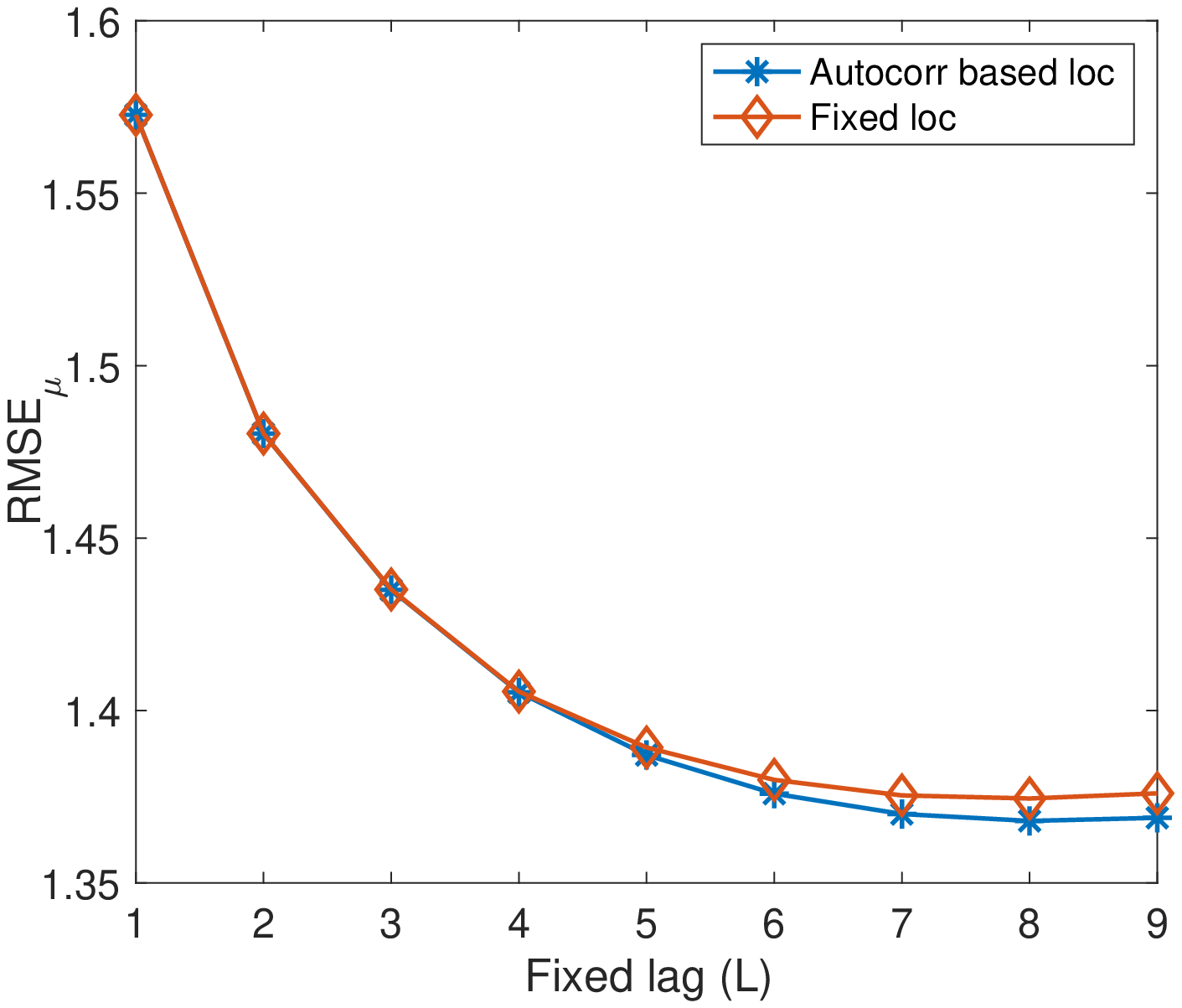} 
		\label{fig:L96_Auto_vs_fix}
	\end{center}
	\caption{Time averaged $\text{RMSE}_{\mu}$ values for ESRS runs with $500000$ assimilation cycles applied to the Lorenz 96 with $M=30$ using two different localisation schemes (fixedloc = stationary scheme and autocorrelation localisation scheme) for localisation radius 1 (left panel) and 8 (right panel).  } 
\end{figure}

%%%%%%%%%%%%%%%%%%%%%%%%%%
The last numerical example we consider is the Lorenz 96 system (L96) \cite{sr:lorenz96} which is given by 
\begin{equation}
\frac{d{\bf x}}{dt}[s]=\Big({\bf x}[s+1]-{\bf x}[s-2]\Big){\bf x}[s-1]-{\bf x}[s]+F
\end{equation}
where $ {\bf x}[{-1}]={\bf x}[{N_x-1}]$, ${\bf x}[0]={\bf x}[{N_x}]$ and  ${\bf x}[N_x+1]={\bf x}[1]$ for each spatial component $s$ and here $N_x=40$. It is another classical toy model which is often utilised to test algorithms in a higher dimensional context. The following setting is used for the Lorenz 96 runs.  We observe every 2nd grid point in space, with temporal observation intervals of $\Delta t_{obs} = 0.11$, observation error variance ${\bf R} = 8$ and forcing constant $F=8$.  A forward Euler scheme is used to discretize the model, with time step $\Delta t = 0.005$.  Additionally, the initial conditions are uncertain, with an ensemble generated from $\mathcal{N}({\bf x}_0, 0.5{\bf I})$.   
Specifically we would like to investigate the spatial as well as the time-wise localisation scheme of the proposed smoothers in this section. The proposed distance functions (\ref{eq:dfix}) and (\ref{eq:dauto}) are compared by means of different ESRS runs where the ensemble size is set to $M=30$ and $500000$ data assimilation steps are computed for each lag ($L=0$ to $L=8$). Further two different localisation radii: ${r}_{loc}=1$ and ${r}_{loc}=8$ are considered to show the effect of a larger radius compared to a very small one. In practice this value is often chosen according to the underlying dynamics of the system.  The distance metric is chosen as $m({\bf p}_{s}, {\bf q}_u):=min\{|s-u|,|s-u-N_x|,|s-u+N_x|\}$, which ensures that the periodic boundary conditions are taken into account. The resulting RMSE values are displayed in the panels of Figure \ref{fig:L96_Auto_vs_fix}. Finally, the Gaspari-Cohn function \cite{sr:gaspari99} is used for $\rho$.

\begin{figure}
\begin{center}
\includegraphics[width=8cm, height=6.857cm]{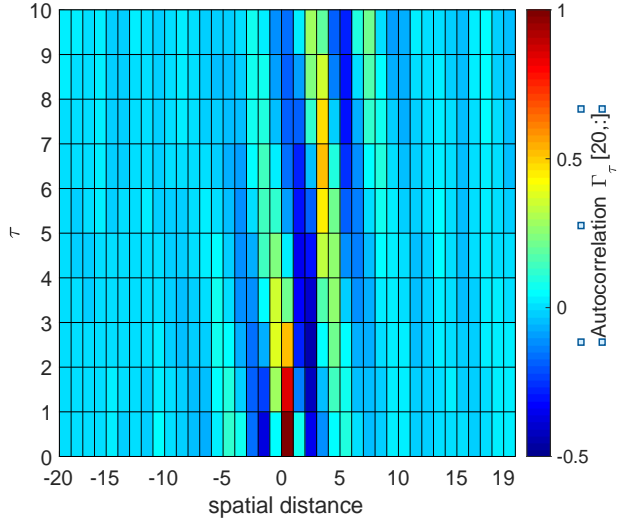} 
\end{center}
\caption{The autocorrelation coefficient $\Gamma_\tau[s,:]$ as defined in (\ref{eq:autocorrcoef}) for the L96 calculated over a time series of 10000 time steps for $s=20$. Distance in space from grid point $s=20$ to all other grid points is shown on the x-axis. Note that one lag unit corresponds to one assimilation step. 
}
\label{fig:autoL96}
\end{figure}

The autocorrelation localisation scheme using the distance function given in (\ref{eq:dauto}) performs equally well as the stationary distance defined in (\ref{eq:dfix}) for lags $0$ to $2$ for both radii due to the fact that $s^*=s$ (see Figure \ref{fig:autoL96}). Yet as soon as $s^*$ and $s$ differ, the autocorrelation localisation scheme produces lower RMSEs than the stationary localisation (see Figure \ref{fig:L96_Auto_vs_fix}). Naturally the choice of radius also effects the performance of the two schemes and the benefits of the localisation scheme based on autocorrelation coefficients is more prominent for radius $1$ as can be seen in the left panel of Figure \ref{fig:L96_Auto_vs_fix}. Thus in particular in a set up with small spatial localisation radii and longer smoother lags, the autocorrelation scheme would be preferred to the stationary one. Yet the stationary scheme is a valid choice for small lags and large spatial localisation radii and does not require the computation of the autocorrelation of the underlying dynamical model.

\begin{figure}
\begin{center}
	\includegraphics[width=8cm, height=6.857cm]{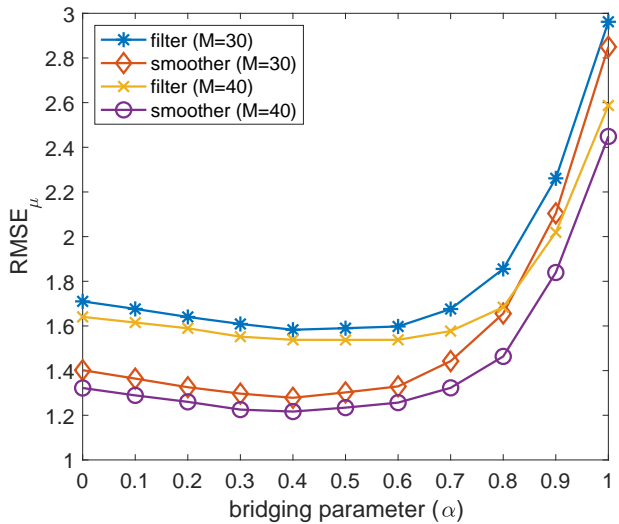} 
\end{center}
\caption{ Time averaged $\text{RMSE}_{\mu}$ values for a range of hybrid parameters $\alpha$ for the ETPF-ESRF (second order corrected) and ETPS-ESRS smoother applied to the L96 model for ensemble sizes $M=30$ and $M=40$. The ETPS corresponds to $\alpha = 1$ and ESRS corresponds to $\alpha = 0$.  Here a stationary localisation scheme was implemented.  } 
\label{fig:RMSE_hybrid_L96}
\end{figure}

Finally a hybrid of an ESRS and the second order corrected ETPS is tested on the L96 system. For the ETPS-ESRS hybrid runs we set $L=10$ and ${r}_{loc}=8$, $\lambda=40$ and consider the stationary localisation scheme based on the distance given in (\ref{eq:dfix}). A total of 50,000 assimilation steps are performed with a burn-in period of $1000$ time instances. Figure \ref{fig:RMSE_hybrid_L96} displays the results for $\alpha\in\{0,0.1,\dots,0.9,1\}$ for $M=30$ and $40$ and the best RMSE values are obtained for $\alpha=0.4$ in the smoothing case. This suggests that a hybrid smoother with localisation could be a useful tool to increase the accuracy of the underlying state estimation problem in a high dimensional setting.

\vspace{2cm}

%%%%%%%%%%%%%%%%%%%%%%%%%%%%%
%
\section{Conclusions}
%
%%%%%%%%%%%%%%%%%%%%%%%%%%%%
A general formulation (\ref{eq:LETS_proper}) for smoothers induced by the family of LETFs is derived. Smoothing extensions for some members of this class were already available e.g., ESRS \cite{jdw:Cosmeetal2009} or NETS \cite{Kirchgessneretal2017}, yet this general extension to a complete class of smoothers allows one to universally investigate the underlying properties of this family rather than having to study them individually.  Furthermore, we show how the LETS framework allows one to develop the extension of the ETPF \cite{acevedo2016} to the smoother case.  Important features for high dimensional state and parameter estimation such as localisation, hybrid implementations and second order corrections are introduced for the LETSs and are studied numerically. 

Application to the Lorenz 63, Mackey-Glass model and Lorenz 96 systems confirmed various improvements that can be achieved via the proposed techniques. Firstly, the proposed second order corrected ETPS has superior performance in terms of RMSEs and CRPS for moderate ensemble sizes.  Similar to previous studies, the ESRS has superior performance for very low ensemble sizes (relative to the state dimension), where it is known that Gaussian assumptions can lead to better tracking capability than attempting to capture the full posterior via importance sampling.  This can be seen in the 40-variable Lorenz96 ($M=30 $ \& $M=40$) and Lorenz63 ($M \leq 20$) experiments.  Secondly, the proposed numerical techniques are shown to lead to improvements on existing smoothers.   
The optimal rotation proposed here significantly improves RMSE scores of the NETS, such that it is almost equivalent to the RMSE of the second order corrected ETPS (see Figure \ref{fig:L63_RMSE_optimallag}) .  In line with earlier results obtained in the context of filtering, the second order corrected ETPS has a significantly higher accuracy than the ETPS without any correction (see Figure \ref{fig:L63_RMSE_optimallag}, right panel).  Additionally, the optimal combination of the ESRS and ETPS through the hybrid framework has better performance in a mean squared error sense than the individual smoothers (see Figures \ref{fig:L63_RMSE_lag} and \ref{fig:RMSE_hybrid_L96}). This suggests that an optimal choice of hybrid smoother could be used to provide more accurate and robust posterior estimates in high dimensional nonlinear applications.

It is crucial to mention however that additional computations associated with smoothing lead to an increase in run time compared to filtering.  It is therefore necessary to find the right trade off between accuracy and computational cost. 
In this regard, we have proposed a number of temporal and spatio-temporal localisation schemes and evaluated them numerically for the Lorenz 96 model (see Figure \ref{fig:L96_Auto_vs_fix}). The results are promising and suggest that smoothing (in particular a hybrid smoother) is a feasible option through appropriate localisation in a high dimensional setting.

%%%%%%%%%%%%%%%%%%%%%%%
%
\section*{Acknowledgement}
%
%%%%%%%%%%%%%%%%%%%%%%%

This research has been funded by 
Deutsche Forschungsgemeinschaft (DFG) through grant 
CRC 1294 \lq\lq Data Assimilation\rq\rq, Project (A02) \lq\lq 
Long-time stability and accuracy of ensemble transform
filter algorithms\rq\rq.

%%%%%%%%%%%%%%%%%%%%%%%%
\bibliographystyle{alpha}
\bibliography{survey_paper}

\end{document}